\documentclass[final,leqno]{siamltex}

\usepackage{times}

\usepackage{amsmath}
\usepackage{amssymb}
\usepackage{algpseudocode}
\usepackage{subfig}
\usepackage{epsfig}
\usepackage[pdftex,letterpaper=true,colorlinks=true,linkcolor=red,filecolor=green,citecolor=red,pdfpagemode=None]{hyperref}
\usepackage{nomencl}
\usepackage{mathrsfs}
\usepackage{chngcntr}
\counterwithout{figure}{section}
\makenomenclature

\newtheorem{example}[theorem]{Example}

\title{On applying the maximum volume principle to a basis selection problem in multivariate polynomial interpolation\thanks{This 
        work was supported by the Academy of Finland (decision 267789).}}

\author{Vesa Kaarnioja\footnotemark[2]}

\begin{document}

\maketitle

\renewcommand{\thefootnote}{\fnsymbol{footnote}}

\footnotetext[2]{Aalto University, Department of Mathematics and Systems Analysis, P.O. Box 11100, FI-00076 Aalto, Finland ({\tt vesa.kaarnioja@aalto.fi}).}

\begin{abstract}
The maximum volume principle is investigated as a means to solve the following problem: Given a set of arbitrary interpolation nodes, how to choose a set of polynomial basis functions for which the Lagrange interpolation problem is well-defined with reasonable interpolation error? The interpolation error is controlled by the Lebesgue constant of multivariate polynomial interpolation and it is proven that the Lebesgue constant can effectively be bounded by the reciprocals of the volume (i.e., determinant in modulus) and the minimal singular value of the multidimensional Vandermonde matrix associated with the interpolation problem. This suggests that a large volume of the Vandermonde system can be used as an indicator of accuracy and stability of the resulting interpolating polynomial. Numerical examples demonstrate that the approach outlined in this paper works remarkably well in practical computations.
\end{abstract}

\begin{keywords} 
Vandermonde matrix, Lebesgue constant, multivariate interpolation, maximum volume principle 
\end{keywords}

\begin{AMS}
41A05, 41A10
\end{AMS}

\pagestyle{myheadings}
\thispagestyle{plain}
\markboth{V. KAARNIOJA}{MAXIMUM VOLUME PRINCIPLE IN POLYNOMIAL BASIS SELECTION}

\section{Introduction}
\label{intro}
The construction of multivariate interpolation rules is usually achieved through tensorization of univariate rules. The complexity of this approach grows exponentially with respect to the dimension of the problem, which makes this approach intractable in problems with even moderate dimensionality. The use of sparse grids~\cite{barthelmann00} reduces this complexity to being polynomial with respect to the dimension. However, node configurations based on tensorized grids are well-defined only as long as the node configuration remains in highly structured format and, in practice, it is difficult to modify the placement of nodes lying on either tensor or sparse grids without compromising the accuracy of the solution.

However, an interpolation node set for multivariate polynomial interpolation need not lie on a tensorized grid to be accurate. The study of optimal node configurations over arbitrary or standard domains---such as the unit disk and square---that produce accurate interpolation formulae goes back to the work of Chung and Yao~\cite{chungyao}. There have been a number of recent developments within this line of research: Sommariva et al.~studied the numerical construction of approximate Fekete points~\cite{briani12,sommarivaqr,sommariva}, Van Barel et al.~studied the construction of nodes obtained by optimization of the Lebesgue constant~\cite{vanbarel14}, and Narayan and Xiu studied the construction of nested nodal sets~\cite{narayanxiuoptimalpoints}. Optimal accuracy nodal points subvert the need to construct costly tensor grids---and thereby the so-called \emph{curse of dimensionality}---without compromising interpolation accuracy, which is a big contributing factor to their appeal in applications such as the polynomial collocation method used for solving parameter-dependent PDEs~\cite{narayanxiu,zitnanart,zitnanproc}.

One reason for the influx of new approximation theory regarding optimal accuracy interpolation nodes stems from the mathematical methods that have recently gained attention in, e.g., data-mining applications, where the CUR matrix approximation is used to obtain low-cost, low-rank approximations of matrices containing immense quantities of data. In particular, the maximum volume principle (i.e., the task of finding the submatrix having maximal determinant in modulus) is an important indicator for finding a quasi-optimal CUR approximation~\cite{maxvol2001}. Several algorithms for the approximate computation of the maximal volume submatrix have been considered in the literature, see, e.g.,~\cite{civrilgreedy} for a description of a greedy algorithm and see~\cite{maxvol} for the \texttt{MaxVol} algorithm.

The maximum volume principle has become a key ingredient in the development of optimal node configurations in multivariate polynomial interpolation: For example, the aforementioned works~\cite{briani12,sommariva,vanbarel14} employ approximate maximum volume Vandermonde submatrices for the identification of nearly optimal accuracy interpolation node configurations. The approach taken in this paper may be regarded as the dual to the problem studied in the aforementioned works: Instead of finding the optimal interpolation nodes with respect to a fixed family of polynomial basis functions, the set of nodes is kept arbitrary and the maximum volume principle is used in the task of finding polynomial basis functions that produce an interpolating polynomial with favorable approximation properties.

\subsection{Related work}

The task of finding an interpolating polynomial for arbitrary node configurations has been considered in the literature by several authors. Kergin interpolation~\cite{micchelli80} provides a constructive, although computationally impractical,  method to develop an interpolating polynomial for any set of nodes. Sauer and Xu~\cite{sauerxu95} proposed an explicit algorithm for incremental addition of points to form the Lagrange and Newton interpolating polynomials under the assumption that the polynomial basis functions are known a priori to produce a well-defined interpolating polynomial; see also~\cite{sauer} for additional developments of this algorithm. De Boor and Ron developed the method of least polynomial interpolation~\cite{deboor1} (see~\cite{deboor2} for computational remarks on this method) which transports the interpolation problem to the dual space of $d$-variate polynomials, i.e., the space of formal $d$-variate Taylor series. Although computationally expensive, this method can be used to produce an interpolating polynomial for arbitrary node configurations. Recently, an extension of least polynomial interpolation was introduced and studied in the framework of the stochastic collocation method by Narayan and Xiu~\cite{narayanxiu}.

More details on the theoretical background and related work concerning multivariate polynomial interpolation can be found in the survey by Gasca and Sauer~\cite{gascasauer00}.

\subsection{Contents of this paper}

This document is organized as follows. The basic notations and preliminaries of multivariate polynomial interpolation are reviewed in Section~\ref{prelims}. The relationship between the Lebesgue constant and the determinant as well as the minimal singular value of the associated Vandermonde system is investigated in Section~\ref{lebesg}, where the main theoretical results of this paper are presented. The methodology of applying the maximum volume principle is examined in numerical experiments in Section~\ref{numex}, and we end with some concluding remarks. Appendix~\ref{sappendix} contains a detailed description of the parameters used in the construction of the Smolyak interpolating polynomial in relation to the numerical example of Subsection~\ref{snodes}.
\section{Notations and preliminaries}\label{prelims}
\subsection{Table of notations}
The special notations used throughout this paper are listed in the following table.

\begin{tabular}{ll}
$\Pi^d$&The space of all real $d$-variate polynomials;\\
$\Pi_k^d$&The space of all real $d$-variate polynomials with total degree at most $k$;\\
$\det_i(M;\mathbf{y})$&The determinant of matrix $M$ with its $i^\text{th}$ column replaced by vector $\mathbf{y}$;\\
$\delta_{i,j}$&The Kronecker symbol, equal to $1$ when $i=j$ and $0$ otherwise;\\
$M_{:,j}$&The column vector corresponding to the $j^\text{th}$ column of matrix $M$;\\
$\mathbf{e}_i$&The $i^\text{th}$ Euclidean standard basis vector;\\
$\sigma_i(M)$&The $i^\text{th}$ singular value of matrix $M$ ordered $\sigma_{i}(M)\geq\sigma_{i+1}(M)$;\\
$\sigma_{\rm max}(M)$&The largest singular value of matrix $M$;\\
$\sigma_{\rm min}(M)$&The smallest singular value of matrix $M$;\\
$\|M\|_2$&The spectral norm of matrix $M$, i.e., $\|M\|_2=\sigma_{\max}(M)$;\\
$\|M\|_{\rm F}$&The Frobenius norm of matrix $M$, i.e., $\|M\|_{\rm F}=(\sum_i\sigma_i(M)^2)^{1/2}$.
\end{tabular}

\nomenclature[1]{$\Pi^d$}{Space of real $d$-variate polynomials}

\makenomenclature

\subsection{Lagrange interpolation problem}
Let $\mathcal{X}=\{\mathbf{x}_1,\ldots,\mathbf{x}_n\}\subset\mathbb{R}^d$ be a set of mutually distinct nodes. The \emph{Lagrange interpolation problem} is to find a polynomial $L_nf\in\Pi^d$ that satisfies
\begin{align}
L_nf(\mathbf{x}_i)=f(\mathbf{x}_i)\quad\text{for }i\in\{1,\ldots,n\}\label{lagrangeproblem}
\end{align}
for any function $f\!:\mathbb{R}^d\to\mathbb{R}$. The problem is well-defined with respect to the polynomial basis $\mathcal{B}=(\phi_i)_{i=1}^n$, where $\phi_i$ are $d$-variate polynomials, if the multidimensional \emph{Vandermonde matrix}
\[
V_{\mathcal{B},\mathcal{X}}=\begin{bmatrix}\phi_1(\mathbf{x}_1)&\cdots&\phi_1(\mathbf{x}_n)\\ \vdots&\ddots&\vdots\\ \phi_n(\mathbf{x}_1)&\cdots&\phi_n(\mathbf{x}_n)\end{bmatrix}
\]
is invertible. Moreover, the solution to~\eqref{lagrangeproblem} can be expressed in terms of the basis functions $\mathcal{B}$ as
\[
L_nf(\mathbf{x})=\sum_{i=1}^nc_i\phi_i(\mathbf{x}),
\]
where the coefficient vector $\mathbf{c}=[c_1,\ldots,c_n]^\textup{T}$ is the solution to the Vandermonde system
\begin{align}
V_{\mathcal{B},\mathcal{X}}^\textup{T} \mathbf{c}=\mathbf{f},\label{vandermondesystem}
\end{align}
where $\mathbf{f}=[f(\mathbf{x}_1),\ldots,f(\mathbf{x}_n)]^\textup{T}$.

Let the elements of the inverse of the Vandermonde matrix be denoted by $w_{i,j}=(V_{\mathcal{B},\mathcal{X}}^{-1})_{i,j}$ for $1\leq i,j\leq n$. Inverting the Vandermonde matrix is equivalent to the determination of the \emph{Lagrange basis}. The Lagrange basis functions can be identified with
\begin{align}
\ell_i^{\mathcal{B},\mathcal{X}}(\mathbf{x})=\mathbf{e}_i^\textup{T} V_{\mathcal{B},\mathcal{X}}^{-1}\begin{bmatrix}\phi_1(\mathbf{x})\\ \vdots\\ \phi_n(\mathbf{x})\end{bmatrix}=\sum_{j=1}^nw_{i,j}\phi_j(\mathbf{x})\label{lagrangerep}
\end{align}
and they satisfy $\ell_i^{\mathcal{B},\mathcal{X}}(\mathbf{x}_j)=\delta_{i,j}$ for $1\leq i,j\leq n$.

Some examples of polynomial bases are given in the following.
\begin{example}
\begin{itemize}
\item[{\rm (i)}] {\rm Tensor products of univariate polynomials $\psi_i$ enumerated by $i\geq 0$ are multivariate polynomials, usually expressed by using multi-index notation $\phi_{\boldsymbol{\alpha}}(x_1,\ldots,x_d)=\psi_{\alpha_1}(x_1)\cdots\psi_{\alpha_d}(x_d)$, where $\boldsymbol{\alpha}=(\alpha_1,\ldots,\alpha_d)\in\mathbb{N}_0^d$. These polynomials can be used to form a basis $\mathcal{B}=(\phi_{\boldsymbol{\alpha}})_{\boldsymbol{\alpha}\in\mathcal{I}}$, where $\mathcal{I}\subset\mathbb{N}_0^d$ is a multisequence subject to, e.g., the degree lexicographic ordering. Fixing an ordering for the multi-indices yields a one-to-one and onto renumbering $\tau\!:\{1,\ldots,\#\mathcal{I}\}\ni i\mapsto\boldsymbol{\alpha}\in \mathcal{I}$ and thus permits enumerating the basis functions as}
\[
\mathcal{B}=(\phi_{\tau(i)})_{i=1}^{\#\mathcal{I}}.
\]
\item[{\rm (ii)}] {\rm The monomial basis for $\Pi_k^2$ is given by}
\[
(\mathbf{x}^{\boldsymbol{\alpha}})_{|\boldsymbol{\alpha}|\leq k}=(1,x_1,x_2,x_1^2,x_1x_2,x_2^2,\ldots,x_1^k,x_1^{k-1}x_2,\ldots,x_1x_2^{k-1},x_2^k),
\]
{\rm where $\boldsymbol{\alpha}\in\mathbb{N}_0^2$, $\mathbf{x}^{\boldsymbol{\alpha}}=x_1^{\alpha_1}x_2^{\alpha_2}$, and the sequence contains $\#\{\boldsymbol{\alpha}\in\mathbb{N}_0^2:|\boldsymbol{\alpha}|\leq k\}=(k^2+3k+2)/2$ basis functions for $k\geq 0$.}
\end{itemize}
\end{example}

Vandermonde matrices are notoriously ill-conditioned and solving the associated system of equations directly is numerically unstable for high values of $n$. However, we make several remarks regarding high-dimensional systems.
\begin{itemize}
\item[(i)] It is the \emph{degree} of the interpolating polynomial that causes ill-conditioning as noted in~\cite{sauer}. For high-dimensional problems, the number of nodes $n$ is related to the polynomial degree $k$ by $n=\binom{k+d}{d}$. In consequence, interpolating polynomials of high degree $k$ for $d\gg 1$ are seldom encountered in practice due to being fundamentally inaccessible from a computational point of view, which mitigates this issue for high-dimensional problems.
\item[(ii)] It is generally preferable to work instead with the better conditioned \emph{Newton basis}. The Vandermonde system can be converted into a Newton system by using LU factorization: Let $V_{\mathcal{B},\mathcal{X}}^\textup{T}=P^{-1}LU$, where left-multiplication by the matrix $P$ is a permutation of rows and $L$ and $U$ are lower and upper triangular matrices, respectively. Then it is sufficient to solve the system
\[
\begin{cases}
L\mathbf{t}=P\mathbf{f},\\
U\mathbf{c}=\mathbf{t}.
\end{cases}
\]
Here, the \emph{transposed} matrix $U^\textup{T}$ describes the change of basis
\[
[\phi_1(\mathbf{x}),\ldots,\phi_n(\mathbf{x})]^\textup{T} = U^\textup{T} [p_1(\mathbf{x}),\ldots,p_n(\mathbf{x})]^\textup{T},
\]
where the \emph{Newton basis functions} $p_i$ satisfy $p_i(\mathbf{\tilde x}_j)=\delta_{i,j}$ for $1\leq j\leq i$, $1\leq i\leq n$, subject to the reordering of the interpolation nodes given by
\[
\begin{bmatrix}
\mathbf{\tilde x}_1^\textup{T}\\ \vdots\\ \mathbf{\tilde x}_n^\textup{T}
\end{bmatrix}=P\begin{bmatrix}
\mathbf{x}_1^\textup{T}\\ \vdots\\ \mathbf{x}_n^\textup{T}\end{bmatrix}.
\]
The above procedure is essentially the matrix counterpart of the incremental Newton interpolation approach described in~\cite{sauer}.
\end{itemize}

Based on remark (i), it is preferable to work with polynomial bases $\mathcal{B}=(\phi_i)_{i=1}^n$ where the ordering of the basis functions is by increasing degree. Moreover, the geometry of the domain and the locations of the interpolation nodes are an important factor in deciding which polynomial basis to use. This paper does not seek to address this question: The burden of the appropriate choice of basis functions for a given node configuration is based on the application and it is left to the user's discretion.

\subsection{Problem setting}
 Let $\mathcal{X}=\{\mathbf{x}_1,\ldots,\mathbf{x}_n\}\subset\mathbb{R}^d$ be a set of mutually distinct nodes. Given a polynomial basis $\mathscr{B}=(\phi_i)_{i=1}^m$ with $m\geq n$, how to determine a subset $\mathcal{B}\subset\mathscr{B}$ for which
 \newcounter{itemtext}
 \begin{itemize}
 \item[(P1)] The Lagrange interpolation problem~\eqref{lagrangeproblem} is well-defined?
 \item[(P2)] The interpolation error has reasonable magnitude (in case there are several bases for which the problem~\eqref{lagrangeproblem} is well-defined)?
 \end{itemize}

 To tackle problem~(P1), let us introduce the generalized Vandermonde matrix $ \mathscr{V}_{\mathscr{B},\mathcal{X}}\in\mathbb{R}^{m\times n}$ defined elementwise by setting 
 \[
 (\mathscr{V}_{\mathscr{B},\mathcal{X}})_{i,j}=\phi_i(\mathbf{x}_j)\quad \text{for }1\leq i\leq m \text{ and }1\leq j\leq n.
 \]
Then the existence of a well-defined Lagrange interpolating polynomial solving~\eqref{lagrangeproblem} is equivalent to finding a Vandermonde submatrix $V_{\mathcal{B},\mathcal{X}}\in\mathbb{R}^{n\times n}$ of $\mathscr{V}_{\mathscr{B},\mathcal{X}}$ such that $\det V_{\mathcal{B},\mathcal{X}}\neq 0$.
 
To find an invertible submatrix for $\mathscr{V}_{\mathscr{B},\mathcal{X}}$, we use the \emph{maximum volume principle}, i.e., we select the $n\times n$ submatrix of $\mathscr{V}_{\mathscr{B},\mathcal{X}}$ which has the largest determinant in modulus out of all possible $n\times n$ submatrices. Of course, while finding the maximum volume submatrix of $\mathscr{V}_{\mathscr{B},\mathcal{X}}$ ensures that we find a basis $\mathcal{B}$ for which the problem~\eqref{lagrangeproblem} is well-defined---provided that such a basis exists in the first place---the problem of finding the actual maximum volume submatrix is NP-hard~\cite{civrilgreedy}. However, the \texttt{MaxVol} algorithm presented in~\cite{maxvol} is a numerically inexpensive way to determine an approximate maximum volume submatrix within supplied tolerance. We refer to the work~\cite{maxvol} for a detailed account on the \texttt{MaxVol} algorithm.

To address problem~(P2), one method to investigate the behavior of the interpolation error is to consider the Lebesgue constant. In the next section, we investigate the relationship between the Lebesgue constant of the interpolating polynomial and the determinant of the associated Vandermonde matrix. In addition, we investigate the connection to the minimum singular value of the Vandermonde matrix.

\section{Bounds on the Lebesgue constant}\label{lebesg}
Let $f\!:\mathbb{R}^d\to\mathbb{R}$ be a continuous function and denote the convex hull of the points $\mathbf{x}_1,\ldots,\mathbf{x}_n\in\mathbb{R}^d$ by
\[
K={\rm conv}(\mathbf{x}_1,\ldots,\mathbf{x}_n)=\left\{\sum_{j=1}^nu_j\mathbf{x}_j:\sum_{j=1}^nu_j=1\text{ and }u_j\geq 0\text{ for all }1\leq j\leq n\right\}.
\]
Let us assume that the Lagrange interpolation problem~\eqref{lagrangeproblem} is well-defined with respect to a polynomial basis $\mathcal{B}=(\phi_i)_{i=1}^n$ that generates the space of polynomials $\mathcal{P}={\rm span}\,\mathcal{B}$ and a set of mutually distinct nodes $\mathcal{X}=\{\mathbf{x}_1,\ldots,\mathbf{x}_n\}\subset\mathbb{R}^d$ such that $\det V_{\mathcal{B},\mathcal{X}}\neq 0$. Then there exists a Lagrange basis $(\ell_i^{\mathcal{B},\mathcal{X}})_{i=1}^n$ such that ${\rm span}(\ell_i^{\mathcal{B},\mathcal{X}})_{i=1}^n=\mathcal{P}$ and the associated interpolation error is controlled by the \emph{Lebesgue constant}
\[
\Lambda_{\mathcal{B},\mathcal{X}}=\max_{\mathbf{x}\in K}\sum_{i=1}^n|\ell_i^{\mathcal{B},\mathcal{X}}(\mathbf{x})|,
\]
which is characterized by the property
\[
\|f-L_nf\|_{L^\infty(K)}\leq (1+\Lambda_{\mathcal{B},\mathcal{X}})\|f-p\|_{L^\infty(K)},
\]
where $p\in\mathcal{P}$ denotes the best polynomial approximation in $\mathcal{P}$ of the continuous function $f$ subject to the uniform norm in $K$. In addition to its application in bounding the error of polynomial interpolation, it is also a measure of the stability of the interpolating polynomial---hence control over the growth of the Lebesgue constant is crucial in high-dimensional problems.

Let us first investigate the relationship between the Lebesgue constant and the minimal singular value of the associated Vandermonde matrix.
\begin{proposition}\label{svbound}
Let $\mathcal{X}=\{\mathbf{x}_1,\ldots,\mathbf{x}_n\}\subset\mathbb{R}^d$ be a mutually distinct set of nodes and $\mathcal{B}=(\phi_i(\mathbf{x}))_{i=1}^n$ a basis of $d$-variate polynomials such that $\det V_{\mathcal{B},\mathcal{X}}\neq 0$. Then the associated Lebesgue constant is bounded by
\[
\Lambda_{\mathcal{B},\mathcal{X}}\leq \frac{Cn}{\sigma_{\min}(V_{\mathcal{B},\mathcal{X}})},
\]
where $C=\max_{\mathbf{x}\in K}\left(\sum_{i=1}^n|\phi_i(\mathbf{x})|^2\right)^{1/2}$.
\end{proposition}
\proof Let us denote the elements of $V_{\mathcal{B},\mathcal{X}}^{-1}$ by $(V_{\mathcal{B},\mathcal{X}}^{-1})_{i,j}=w_{i,j}$ for $1\leq i,j\leq n$. Recalling the identity~\eqref{lagrangerep}, we can use the Cauchy--Schwarz inequality to obtain
\[
|\ell_i^{\mathcal{B},\mathcal{X}}(\mathbf{x})|=\left|\sum_{j=1}^n w_{i,j}\phi_j(\mathbf{x})\right|\leq \left(\sum_{j=1}^n|w_{i,j}|^2\right)^{1/2}\!\left(\sum_{j=1}^n|\phi_j(\mathbf{x})|^2\right)^{1/2}\leq C\left(\sum_{j=1}^n|w_{i,j}|^2\right)^{1/2}\!.
\]
A second application of the Cauchy--Schwarz inequality yields
\[
\Lambda_{\mathcal{B},\mathcal{X}}=\max_{\mathbf{x}\in K}\sum_{i=1}^n|\ell_i^{\mathcal{B},\mathcal{X}}(\mathbf{x})|\leq C\sqrt{n}\left(\sum_{i,j=1}^n|w_{i,j}|^2\right)^{1/2}=C\sqrt{n}\|V_{\mathcal{B},\mathcal{X}}^{-1}\|_{\rm F}.
\]
The claim follows by utilizing the equivalence of the matrix norms $\|M\|_2\leq \|M\|_{\rm F}\leq \sqrt{n}\|M\|_2$, the identity $\|M\|_2=\sigma_{\max}(M)$, and the fact that $\|V_{\mathcal{B},\mathcal{X}}^{-1}\|_2=1/\sigma_{\min}(V_{\mathcal{B},\mathcal{X}})$.\quad \endproof

A consequence of Proposition~\ref{svbound} is that choosing a submatrix $V_{\mathcal{B},\mathcal{X}}$ of $\mathscr{V}_{\mathscr{B},\mathcal{X}}$ which has the largest minimal singular value (\texttt{MaxMinSv}) out of all $\mathcal{B}\subset\mathscr{B}$ is a good candidate in terms of ensuring reasonable interpolation accuracy and stability of the interpolating polynomial. Unfortunately, there do not appear to exist any efficient algorithms in the literature that are designed to find the approximate \texttt{MaxMinSv} submatrix \emph{directly}. However, something can be said on the \emph{indirect} approximation of the \texttt{MaxMinSv} submatrix.

\begin{itemize}
\item[(i)] It was noted in~\cite{maxminsv} that the problem of determining the \texttt{MaxMinSv} submatrix can be replaced with the problem of determining the \texttt{MaxVol} submatrix, since both kinds of submatrices have qualitatively similar approximation properties.
\item[(ii)] The remark (i) is backed up by intuition: The modulus of the determinant is the product of all singular values. By seeking out the \texttt{MaxVol} submatrix, we do not expect the minimal singular value of such a matrix to be unreasonably small either.
\end{itemize}

The Lebesgue constant can also be related to the determinant of the maximum volume submatrix as the following corollary shows.
\begin{corollary}
Under the assumptions of Proposition~\ref{svbound}, we have
\[
\Lambda_{\mathcal{B},\mathcal{X}}\leq \frac{CD\sqrt{\rm e}\,n}{|\det V_{\mathcal{B},\mathcal{X}}|},
\]
where ${\rm e}$ denotes Euler's number and
\[
D=\max\left\{\sqrt{\frac{\min_{1\leq i\leq n}\sum_{j=1}^n|\phi_j(\mathbf{x}_i)|^2}{\prod_{i=1}^n\sum_{j=1}^n|\phi_j(\mathbf{x}_i)|^2}},\sqrt{\frac{\min_{1\leq i\leq n}\sum_{j=1}^n|\phi_i(\mathbf{x}_j)|^2}{\prod_{i=1}^n\sum_{j=1}^n|\phi_i(\mathbf{x}_j)|^2}}\right\}^{-1}.
\]
\end{corollary}
\proof By~\cite[Theorem~1]{lowerbound}, the minimal singular value can be bounded from below by
\[
\sigma_{\min}(V_\mathcal{B})\geq \left(\frac{n-1}{n}\right)^{(n-1)/2}|\det V_\mathcal{B}|\, D^{-1}.
\]
The assertion follows by plugging this expression into Proposition~\ref{svbound} and observing that $\left(\frac{n}{n-1}\right)^{(n-1)/2}\leq\sqrt{\rm e}$ for all $n\geq 1$.\quad\endproof

An intriguing question encountered in practice is the effect that the addition of a new node and basis function to an existing interpolating polynomial has on the Lebesgue constant. We investigate this question under the assumption that the new node lies in the convex hull of the existing interpolation node configuration.

\begin{proposition}
Let $\mathcal{X}=\{\mathbf{x}_1,\ldots,\mathbf{x}_n\}\subset\mathbb{R}^d$ be  mutually distinct nodes  and $\mathcal{B}=(\phi_i)_{i=1}^n$ a basis of $d$-variate polynomials for which $\det V_{\mathcal{B},\mathcal{X}}\neq 0$. Let $\mathbf{x}_{n+1}\in {\rm conv}(\mathbf{x}_1,\ldots,\mathbf{x}_n)=K$ be a node, and $\phi_{n+1}\in\Pi^d$ a basis function such that $\mathbf{x}_{n+1}\neq \mathbf{x}_i$ and $\phi_{n+1}\neq\phi_i$ for all $i\in\{1,\ldots,n\}$. Let $\det V_{\mathcal{B}\cup\{\mathbf{x}_{n+1}\},\mathcal{X}\cup\{\mathbf{x}_{n+1}\}}\neq 0$. Then
\[
\Lambda_{\mathcal{B}\cup\{\phi_{n+1}\},\mathcal{X}\cup\{\mathbf{x}_{n+1}\}}\leq \Lambda_{\mathcal{B},\mathcal{X}}+\|\ell_{n+1}^{\mathcal{B}\cup\{\phi_{n+1}\},\mathcal{X}\cup\{\mathbf{x}_{n+1}\}}\|_{L^\infty(K)}(1+\Lambda_{\mathcal{B},\mathcal{X}}),
\]
where
\[
\|\ell_{n+1}^{\mathcal{B}\cup\{\phi_{n+1}\},\mathcal{X}\cup\{\mathbf{x}_{n+1}\}}\|_{L^\infty(K)}\leq\left|\frac{\det V_{\mathcal{B},\mathcal{X}}}{\det V_{\mathcal{B}\cup\{\phi_{n+1}\},\mathcal{X}\cup\{\mathbf{x}_{n+1}\}}}\right|\|\phi_{n+1}\|_{L^\infty(K)}(1+\Lambda_{\mathcal{B},\mathcal{X}}).
\]
\end{proposition}
\proof For ease of presentation, let us denote $V_n=V_{\mathcal{B},\mathcal{X}}$, $V_{n+1}=V_{\mathcal{B}\cup\{\phi_{n+1}\},\mathcal{X}\cup\{\mathbf{x}_{n+1}\}}$, $\ell_i^n=\ell_i^{\mathcal{B},\mathcal{X}}$, and $\ell_i^{n+1}=\ell_i^{\mathcal{B}\cup\{\phi_{n+1}\},\mathcal{X}\cup\{\mathbf{x}_{n+1}\}}$, respectively, and let $\mathbf{x}\in K$. The updated Lagrange basis functions are given constructively by the formulae
\begin{align*}
\ell_{n+1}^{n+1}(\mathbf{x})=\displaystyle\frac{\phi_{n+1}(\mathbf{x})-\sum_{i=1}^n\ell_i^n(\mathbf{x})\phi_{n+1}(\mathbf{x}_i)}{\phi_{n+1}(\mathbf{x}_{n+1})-\sum_{i=1}^n\ell_i^n(\mathbf{x}_{n+1})\phi_{n+1}(\mathbf{x}_i)},\\
\ell_i^{n+1}(\mathbf{x})=\ell_i^n(\mathbf{x})-\ell_{n+1}^{n+1}(\mathbf{x})\ell_i^n(\mathbf{x}_{n+1})\quad\text{for }1\leq i\leq n.
\end{align*}
It follows that
\begin{align*}
\sum_{i=1}^{n+1}|\ell_i^{n+1}(\mathbf{x})|&=\sum_{i=1}^n|\ell_i^n(\mathbf{x})-\ell_{n+1}^{n+1}(\mathbf{x})\ell_i^n(\mathbf{x}_{n+1})|+|\ell_{n+1}^{n+1}(\mathbf{x})|\\
&\leq\sum_{i=1}^n|\ell_i^n(\mathbf{x})|+\|\ell_{n+1}^{n+1}\|_{L^\infty(K)}\sum_{i=1}^n|\ell_i^n(\mathbf{x}_{n+1})|+\|\ell_{n+1}^{n+1}\|_{L^\infty(K)}\\
&\leq \Lambda_{\mathcal{B},\mathcal{X}}+\|\ell_{n+1}^{n+1}\|_{L^\infty(K)}(1+\Lambda_{\mathcal{B},\mathcal{X}}),
\end{align*}
where the upper bound on $\Lambda_{\mathcal{B}\cup\{\phi_{n+1}\},\mathcal{X}\cup\{\mathbf{x}_{n+1}\}}$ follows by taking the maximum over $\mathbf{x}\in K$.

On the other hand, Schur's determinant identity yields
\begin{align*}
{\rm det}\sb{n+1}(V_{n+1};\mathbf{x})&=\det\!\left[\begin{array}{c|c}
V_n&\begin{array}{c}\phi_1(\mathbf{x})\\ \vdots\\ \phi_n(\mathbf{x})\end{array}\\
\hline\begin{array}{ccc}\phi_{n+1}(\mathbf{x}_1)&\cdots&\phi_{n+1}(\mathbf{x}_n)\end{array}&\phi_{n+1}(\mathbf{x})
\end{array}\right]\\
&=\det V_n\!\left(\phi_{n+1}(\mathbf{x})-[\phi_{n+1}(\mathbf{x}_1),\ldots,\phi_{n+1}(\mathbf{x}_n)]V_n^{-1}\begin{bmatrix}\phi_1(\mathbf{x})\\ \vdots\\ \phi_n(\mathbf{x})\end{bmatrix}\right)\\
&=\det V_n\!\left(\phi_{n+1}(\mathbf{x})-\sum_{i=1}^n\ell_i^n(\mathbf{x})\phi_{n+1}(\mathbf{x}_i)\right),
\end{align*}
where the final equality follows from~\eqref{lagrangerep}. By taking absolute values on both sides we obtain
\begin{align*}
|{\rm det}\sb{n+1}(V_{n+1};\mathbf{x})|&=|\det V_n|\bigg|\phi_{n+1}(\mathbf{x})-\sum_{i=1}^n\ell_i^n(\mathbf{x})\phi_{n+1}(\mathbf{x}_i)\bigg|\\
&\leq|\det V_n|\|\phi_{n+1}\|_{L^\infty(K)}(1+\Lambda_{\mathcal{B},\mathcal{X}})
\end{align*}
and by utilizing Cramer's rule we have $\ell_{n+1}^{n+1}(\mathbf{x})\det V_{n+1}=\det_{n+1}(V_{n+1};\mathbf{x})$, thus
\[
\|\ell_{n+1}^{n+1}\|_{L^\infty(K)}\leq\left|\frac{\det V_n}{\det V_{n+1}}\right|\|\phi_{n+1}\|_{L^\infty(K)}(1+\Lambda_{\mathcal{B},\mathcal{X}})
\]
proving the assertion.\quad\endproof

\section{Numerical experiments}\label{numex}

We consider the Lebesgue constant of both random and deterministic node configurations in the numerical experiments. The dimensionality of the experiments is kept low with $d\in\{2,3\}$ since this enables reliable discretization of the convex hull generated by the node configurations.

The Lebesgue constant is approximated by creating a triangular mesh for the convex hull $K={\rm conv}(\mathbf{x}_1,\ldots,\mathbf{x}_n)$ of the interpolation nodes $\mathcal{X}=\{\mathbf{x}_1,\ldots,\mathbf{x}_n\}\subset\mathbb{R}^d$ and computing the discrete approximation
\[
\Lambda_{\mathcal{B},\mathcal{X}}\approx \max_{1\leq j\leq N}\sum_{i=1}^n|\ell_i^{\mathcal{B},\mathcal{X}}(\mathbf{y}_j)|
\]
over the mesh vertices $\{\mathbf{y}_1,\ldots,\mathbf{y}_N\}\subset K$ for a given polynomial basis $\mathcal{B}$. The Lagrange basis functions are computed by solving the Vandermonde system~\eqref{vandermondesystem} for $\mathbf{f}\in\{\mathbf{e}_1,\ldots,\mathbf{e}_n\}$ and utilizing formula~\eqref{lagrangerep}.

In Subsection~\ref{randomnodes}, each realization of a sequence of random nodes determines a unique convex hull, which is discretized into a triangular mesh when $d=2$ and a tetrahedral mesh when $d=3$ with each cell not exceeding $10^{-2}$ in Lebesgue measure. In Subsection~\ref{snodes}, the convex hull is the set $[-1,1]^d$ for $d\in\{2,3\}$, which is likewise discretized into a triangular or tetrahedral mesh with each cell not exceeding $10^{-2}$ in Lebesgue measure.

\subsection{Random nodes}\label{randomnodes}

\begin{figure}[!h]
\captionsetup[subfigure]{labelformat=empty}
\centering
\subfloat{{\includegraphics[width=.33 \textwidth]{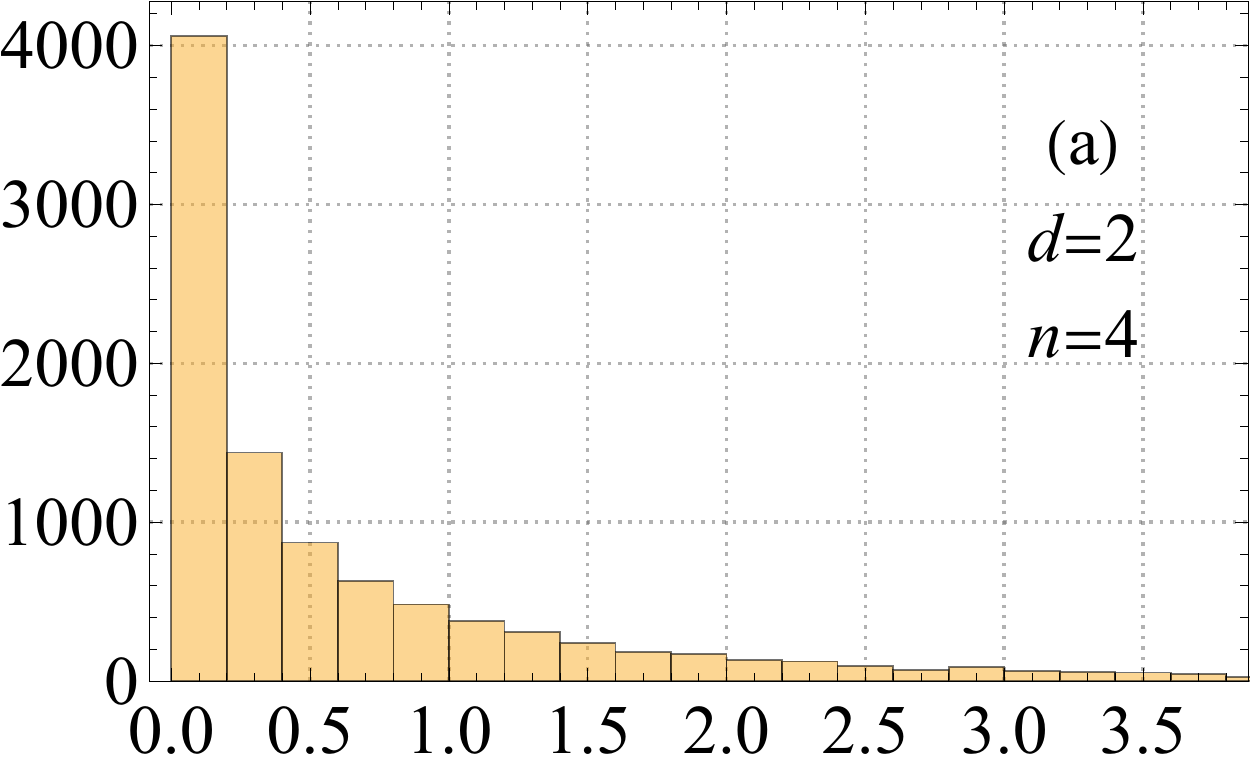}}}
\subfloat{{\includegraphics[width=.33 \textwidth]{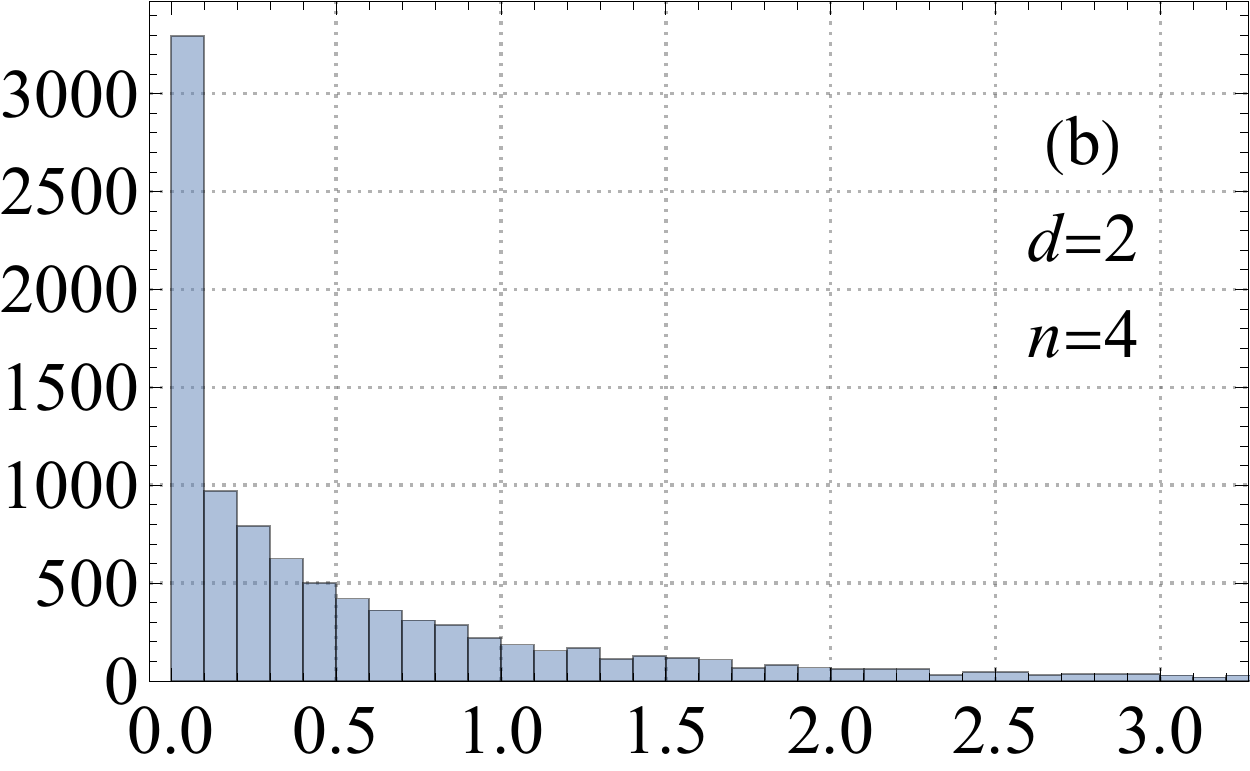}}}
\subfloat{{\includegraphics[width=.33 \textwidth]{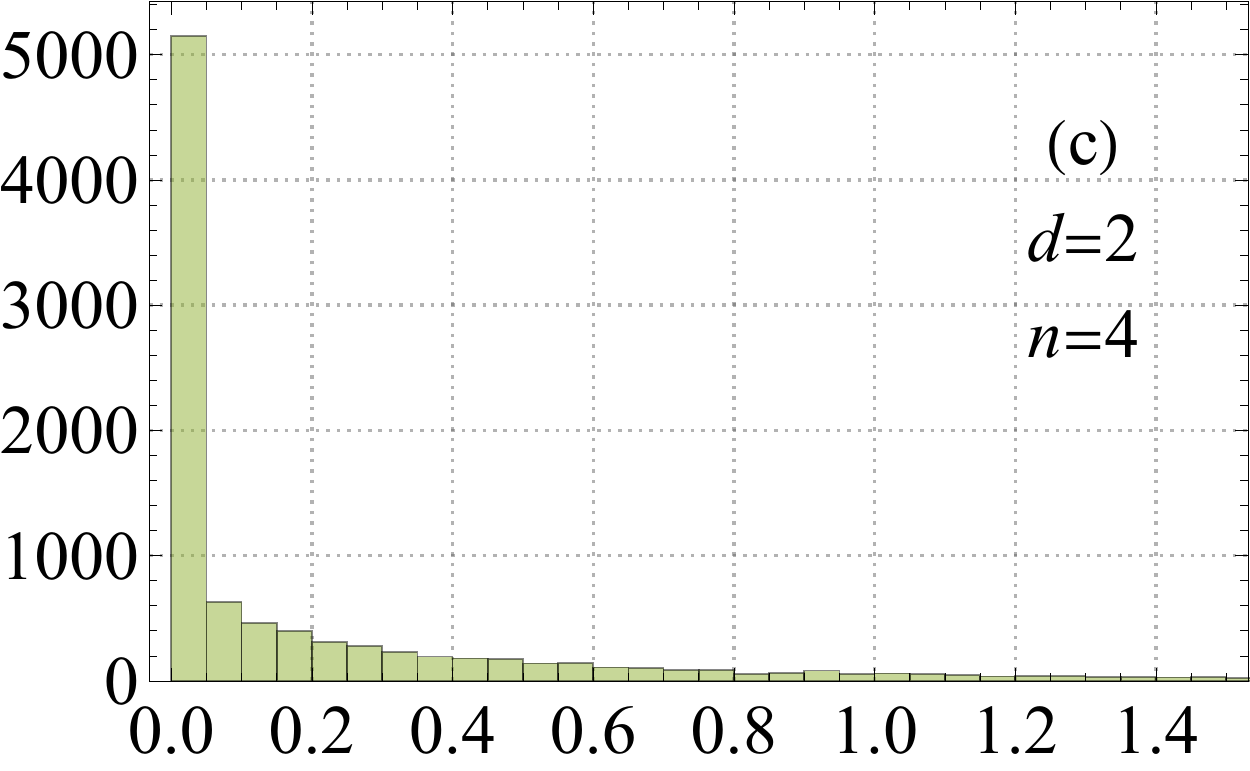}}}\\
\subfloat[$|\Lambda_{d,n}^{\rm best}-\Lambda_{d,n}^{\rm MaxVol}|$]{{\includegraphics[width=.33 \textwidth]{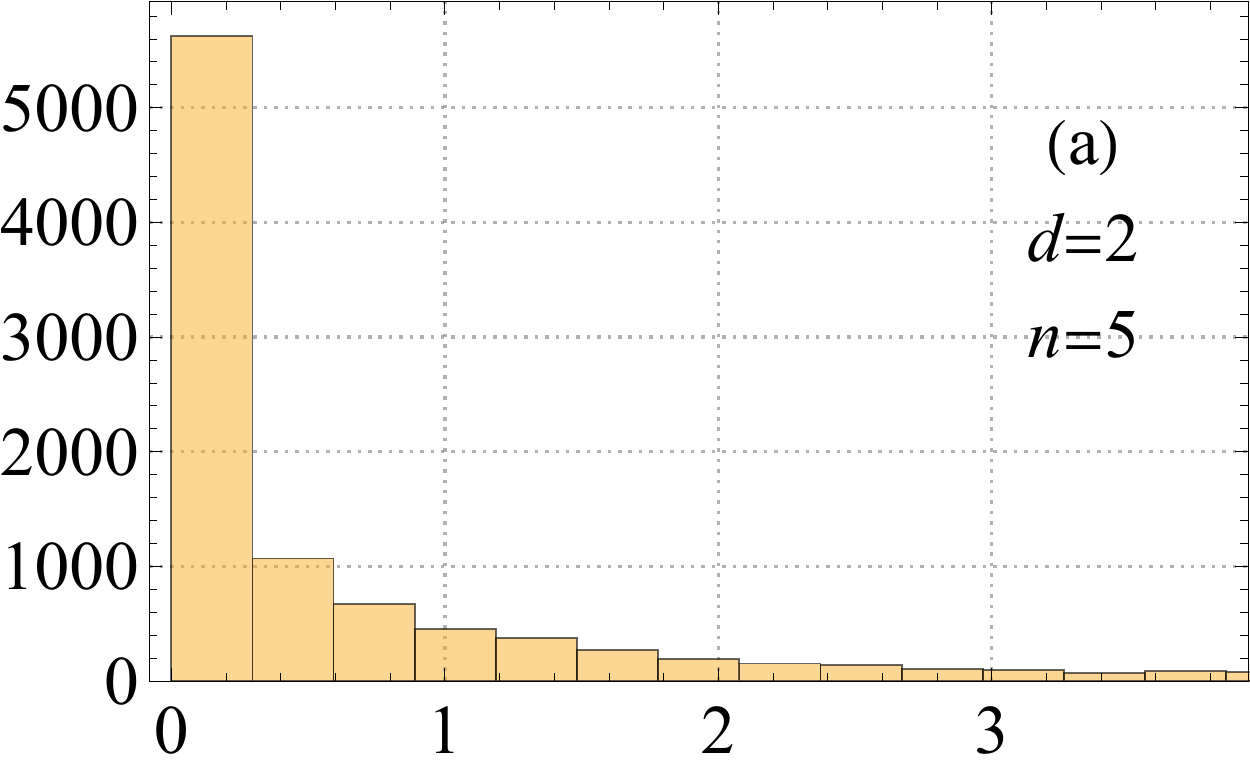}}}
\subfloat[$|\Lambda_{d,n}^{\rm best}-\Lambda_{d,n}^{\rm MaxMinSv}|$]{{\includegraphics[width=.33 \textwidth]{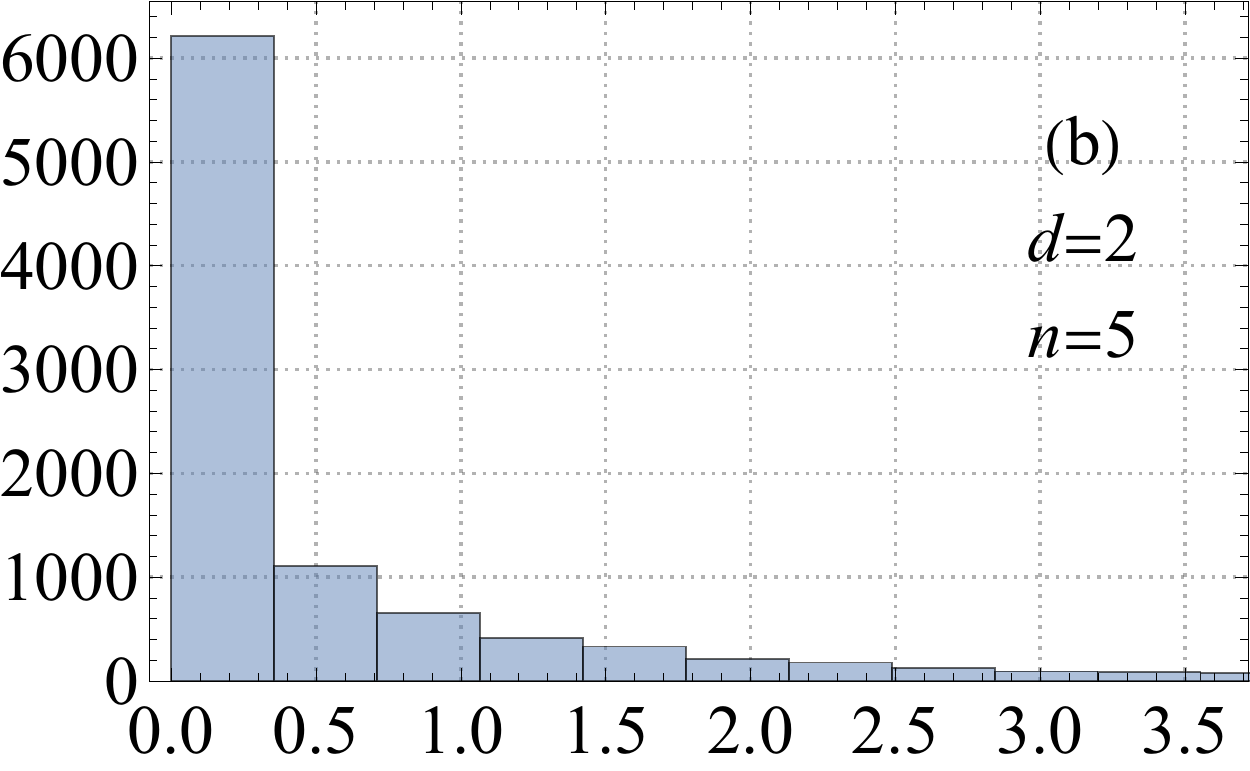}}}
\subfloat[$|\Lambda_{d,n}^{\rm MaxVol}-\Lambda_{d,n}^{\rm MaxMinSv}|$]{{\includegraphics[width=.33 \textwidth]{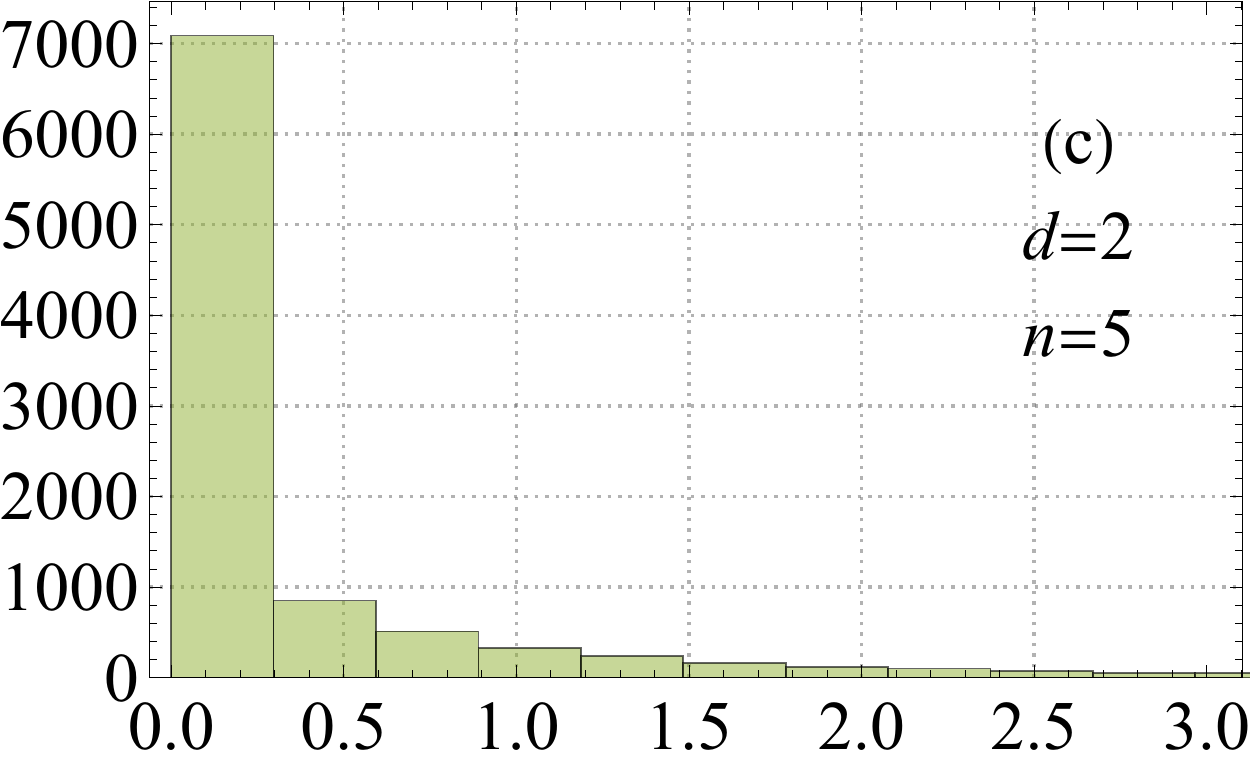}}}
\caption{Case (i), where $10\,000$ realizations of uniformly random node sequences with $d=2$ were generated for each $n\in\{4,5\}$.  The occurrences of the differences (a)--(c) between the Lebesgue constants corresponding to differently chosen polynomial bases have been tabulated in the histograms above.}\label{testi}
\end{figure}

\begin{figure}[!h]
\captionsetup[subfigure]{labelformat=empty}
\centering
\subfloat{{\includegraphics[width=.33 \textwidth]{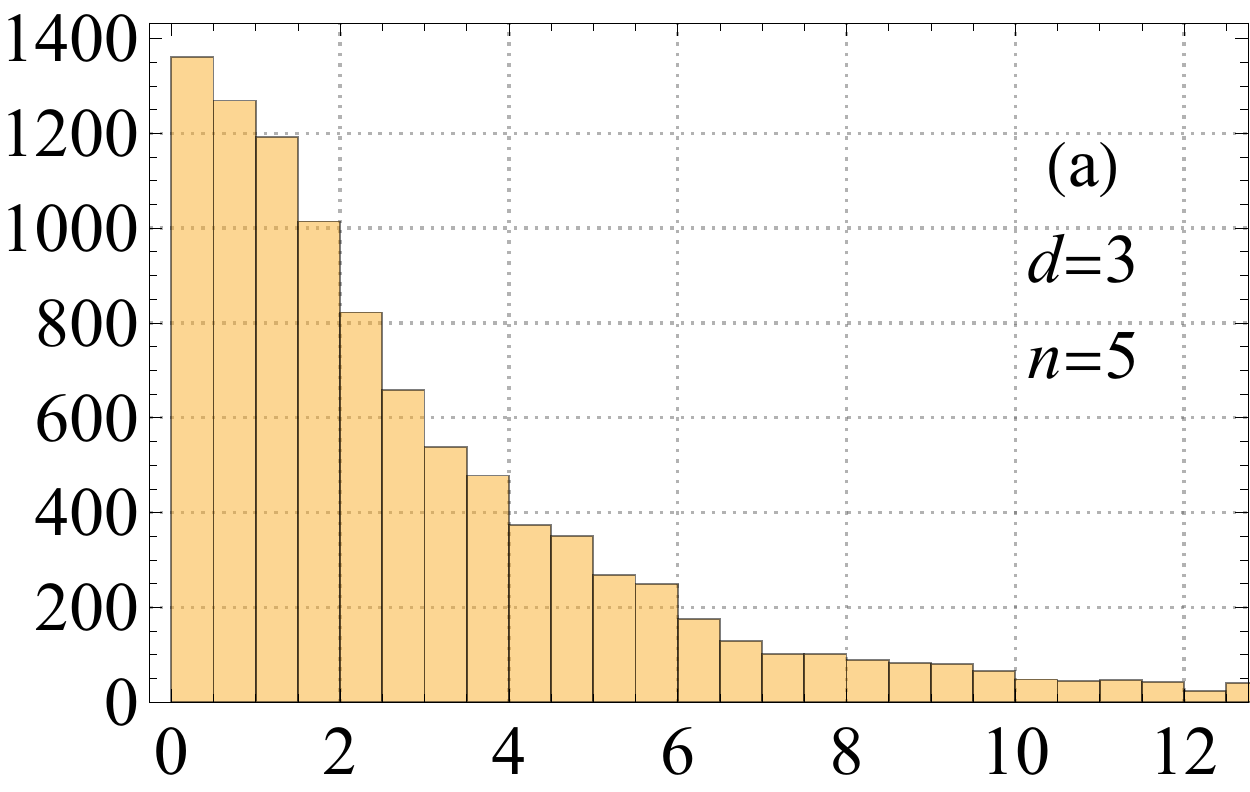}}}
\subfloat{{\includegraphics[width=.33 \textwidth]{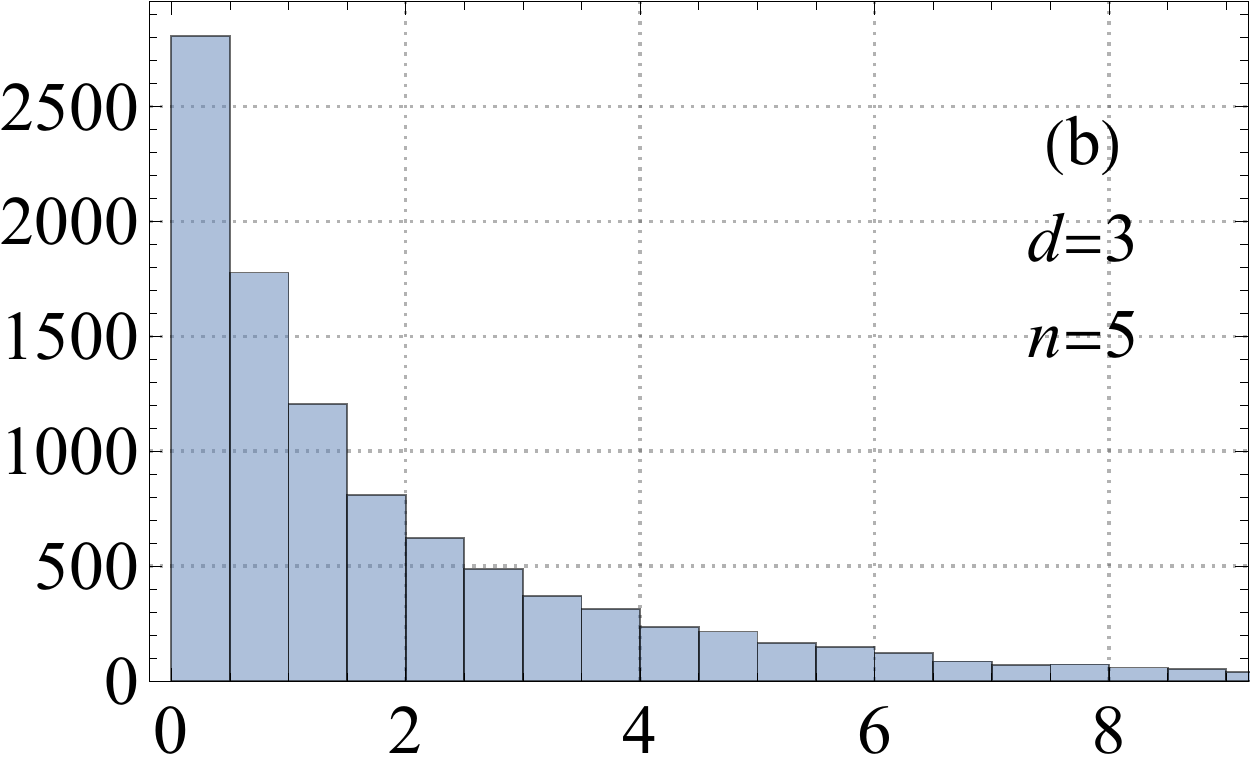}}}
\subfloat{{\includegraphics[width=.33 \textwidth]{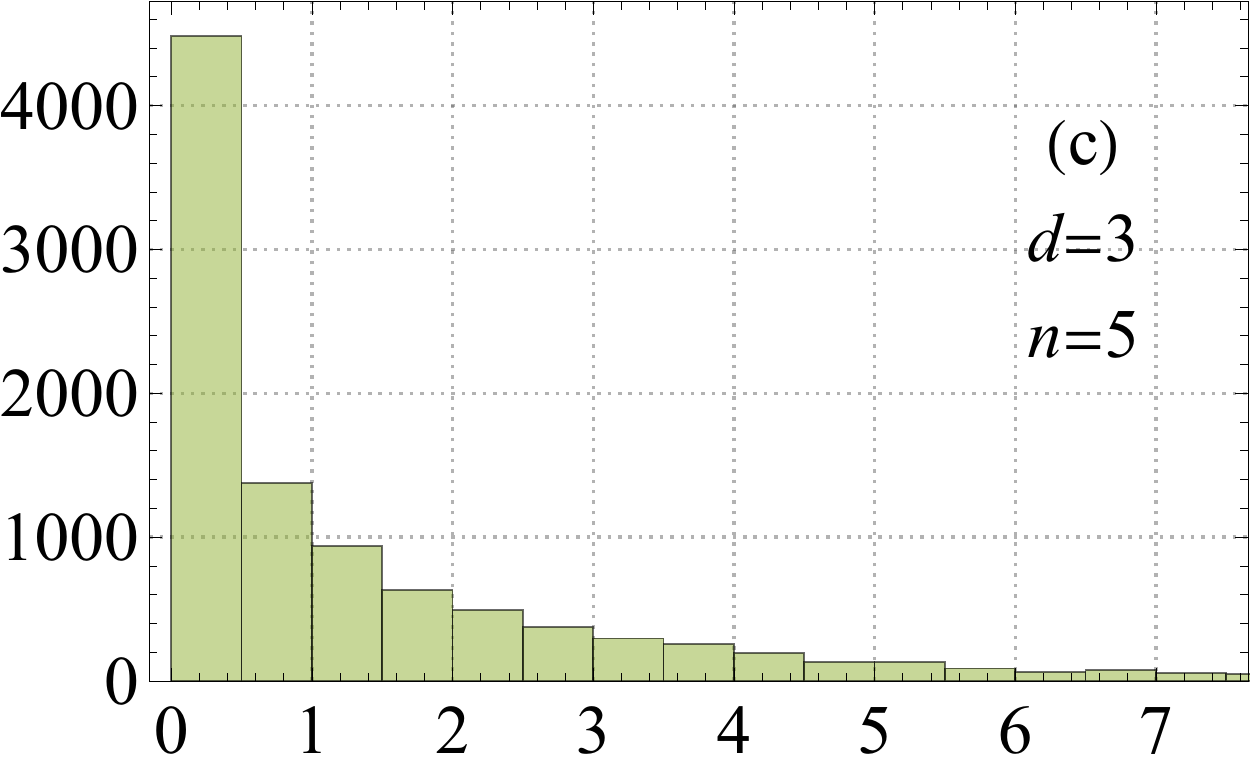}}}\\
\subfloat{{\includegraphics[width=.33 \textwidth]{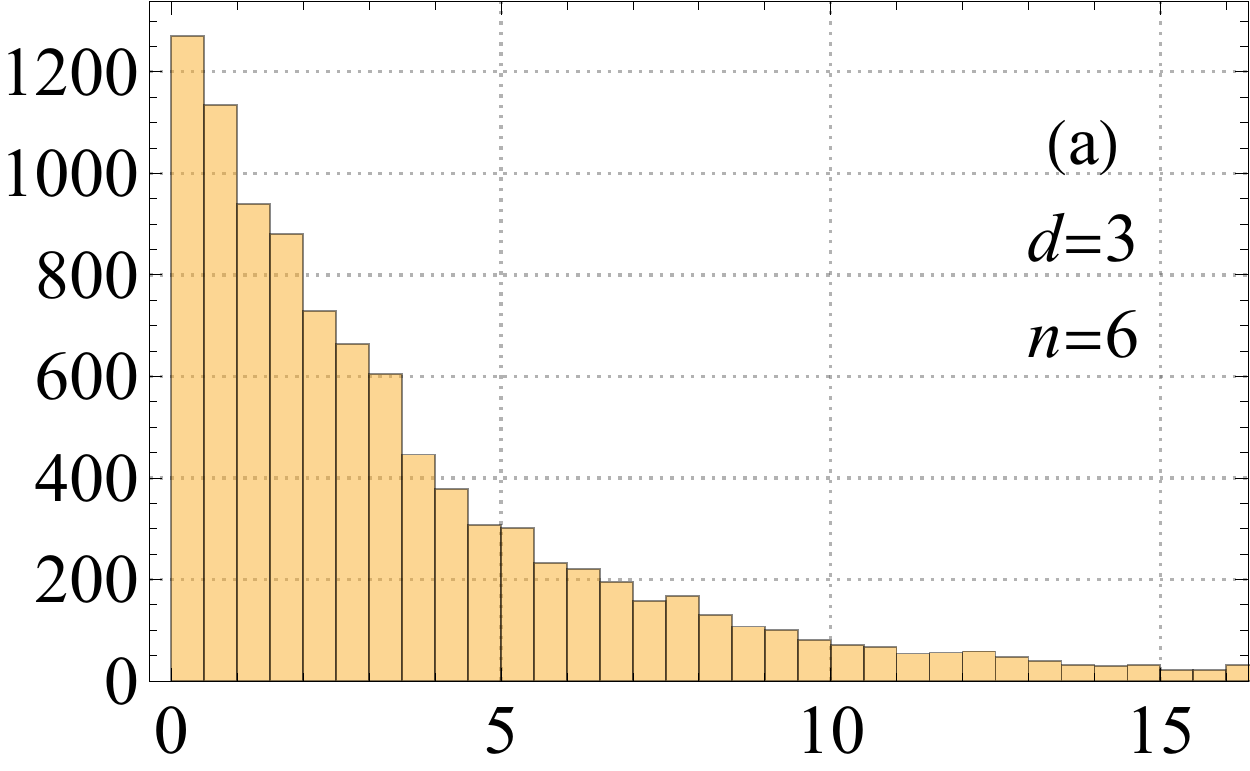}}}
\subfloat{{\includegraphics[width=.33 \textwidth]{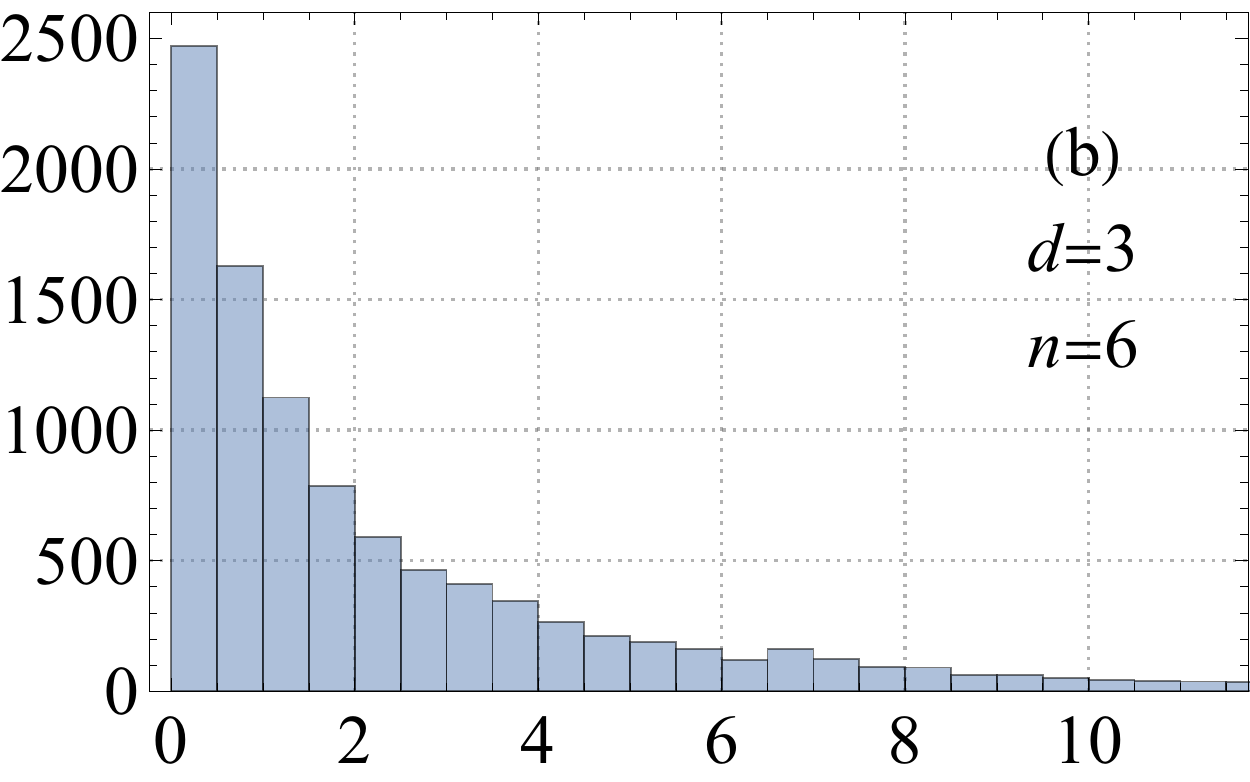}}}
\subfloat{{\includegraphics[width=.33 \textwidth]{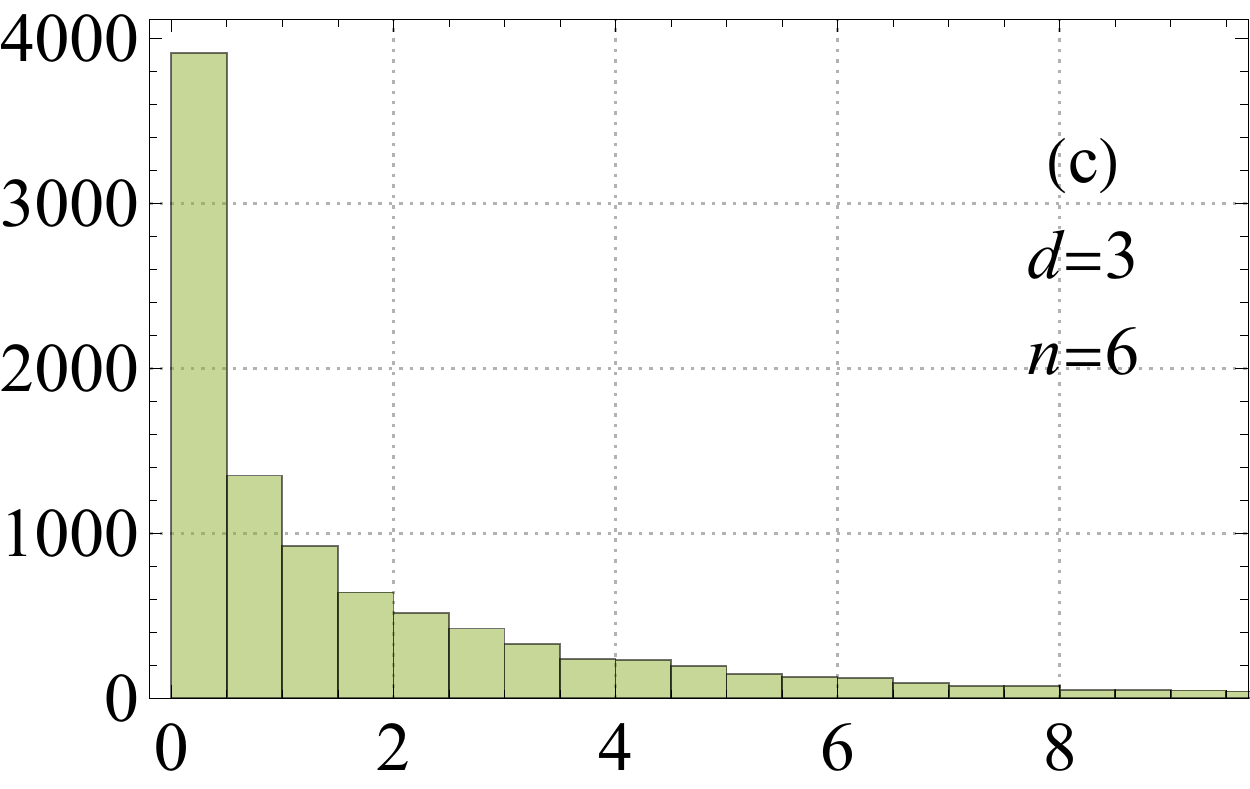}}}\\
\subfloat{{\includegraphics[width=.33 \textwidth]{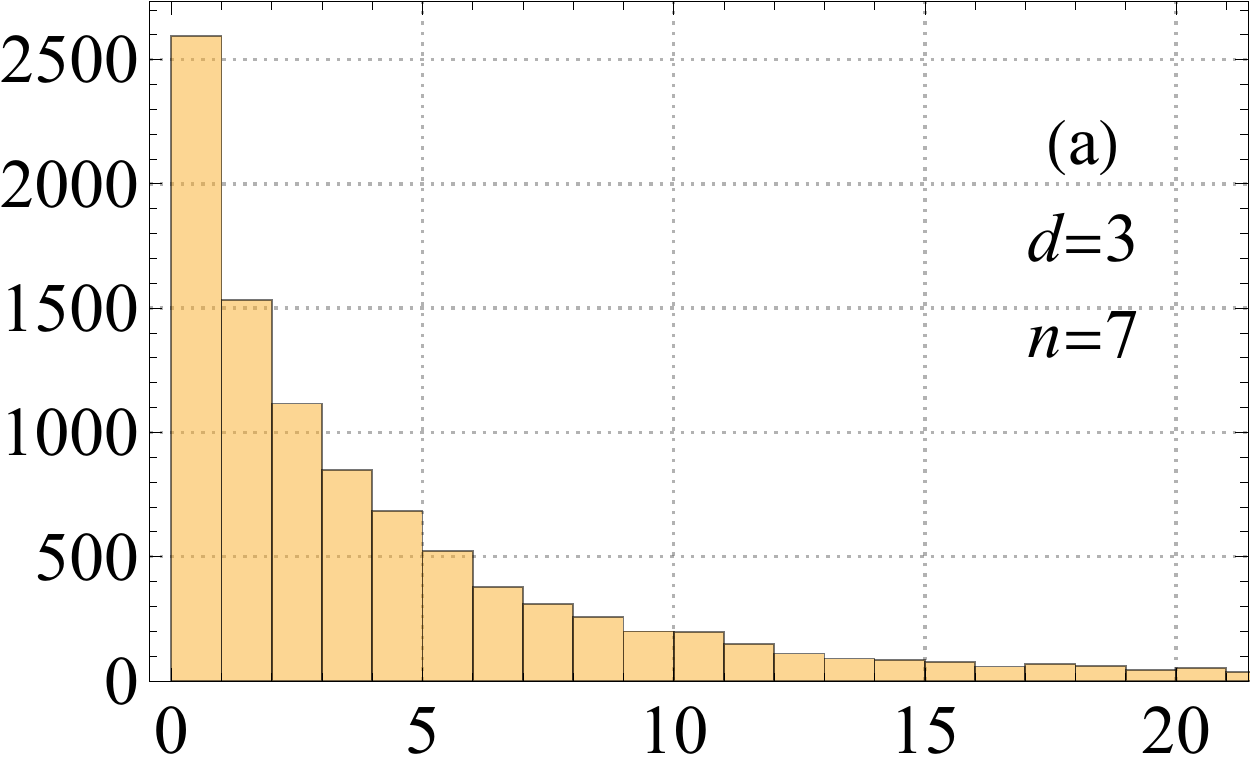}}}
\subfloat{{\includegraphics[width=.33 \textwidth]{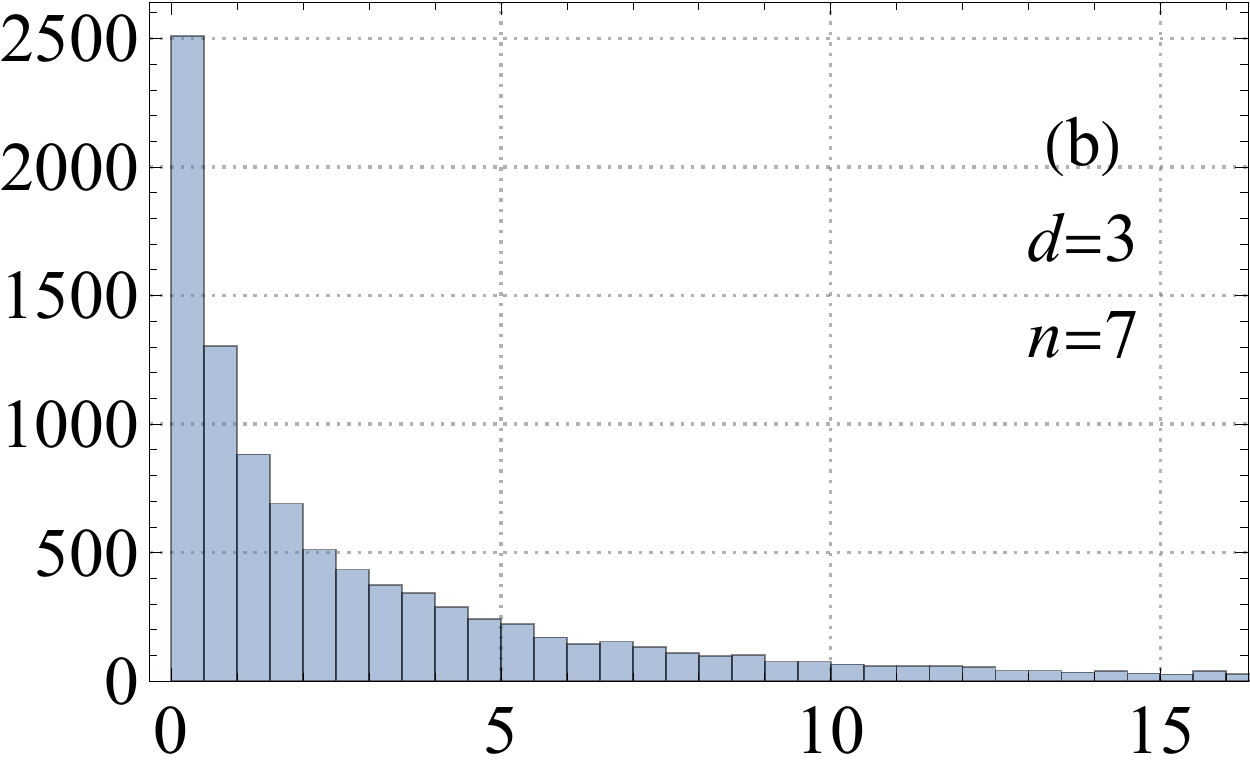}}}
\subfloat{{\includegraphics[width=.33 \textwidth]{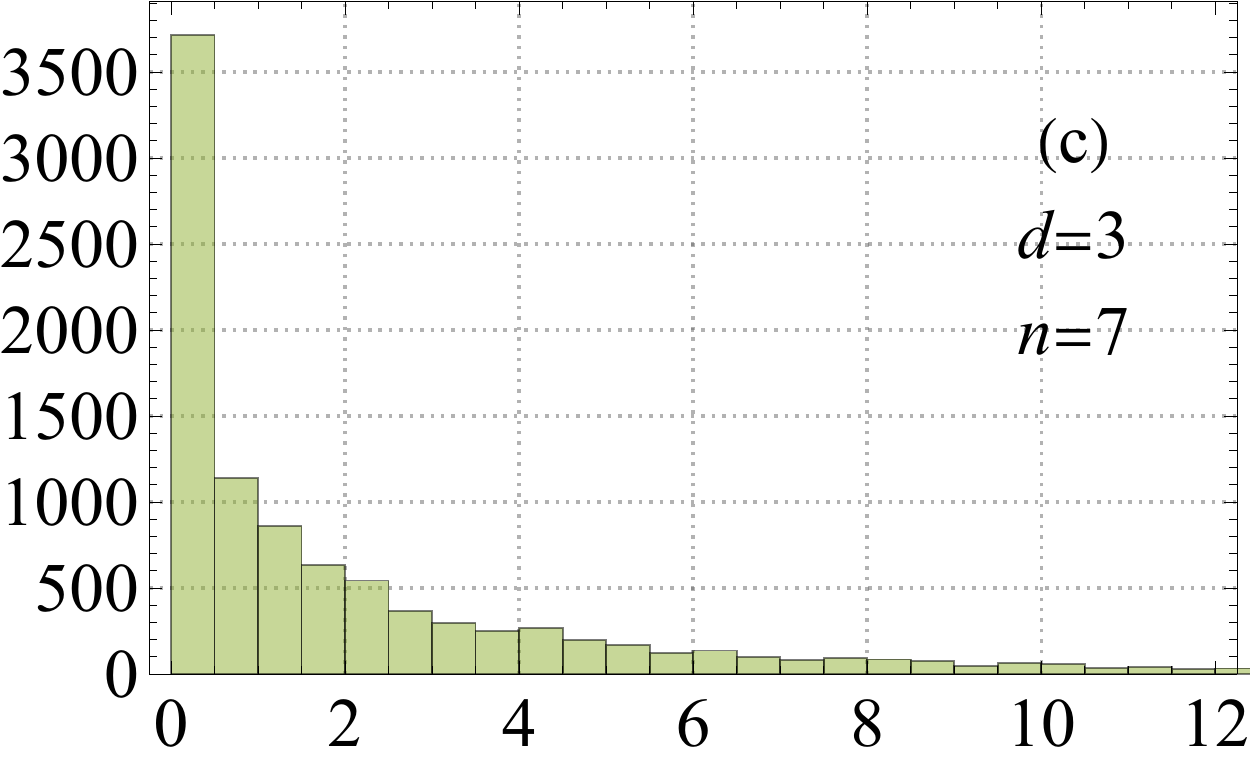}}}\\
\subfloat{{\includegraphics[width=.33 \textwidth]{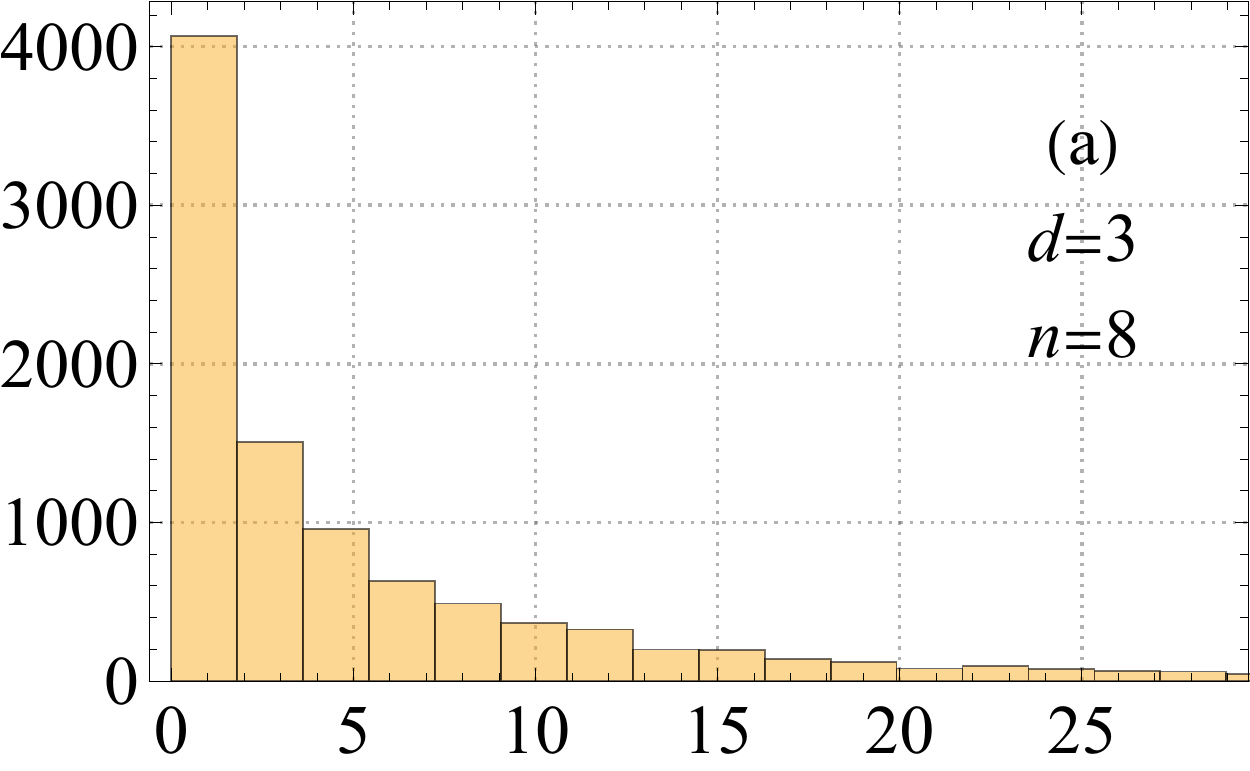}}}
\subfloat{{\includegraphics[width=.33 \textwidth]{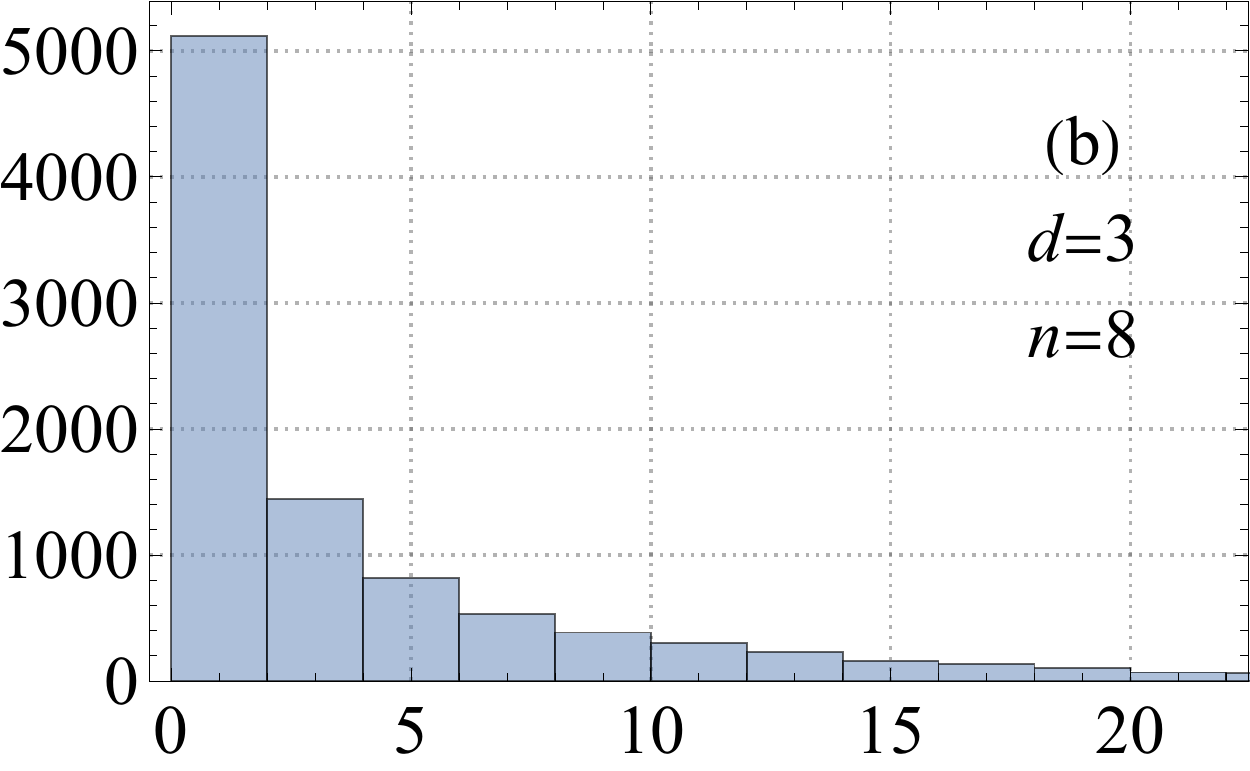}}}
\subfloat{{\includegraphics[width=.33 \textwidth]{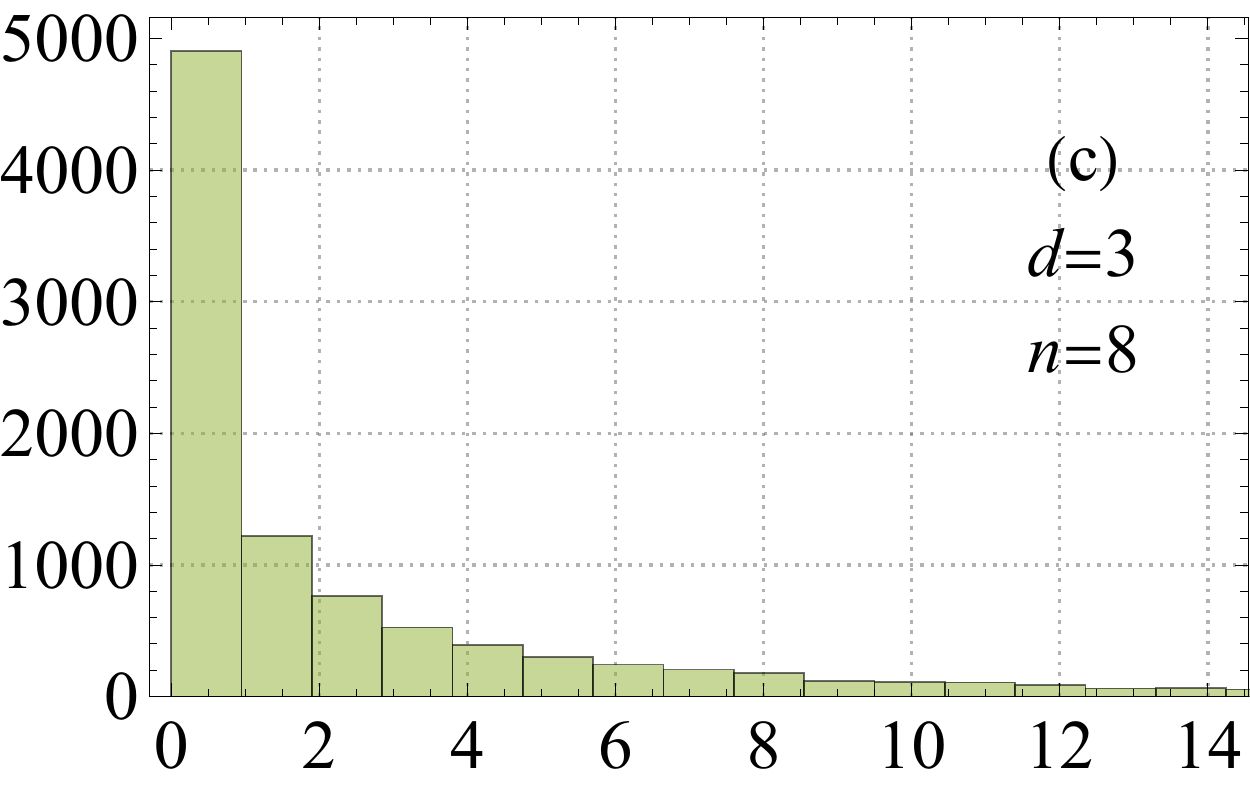}}}\\
\subfloat[$|\Lambda_{d,n}^{\rm best}-\Lambda_{d,n}^{\rm MaxVol}|$]{{\includegraphics[width=.33 \textwidth]{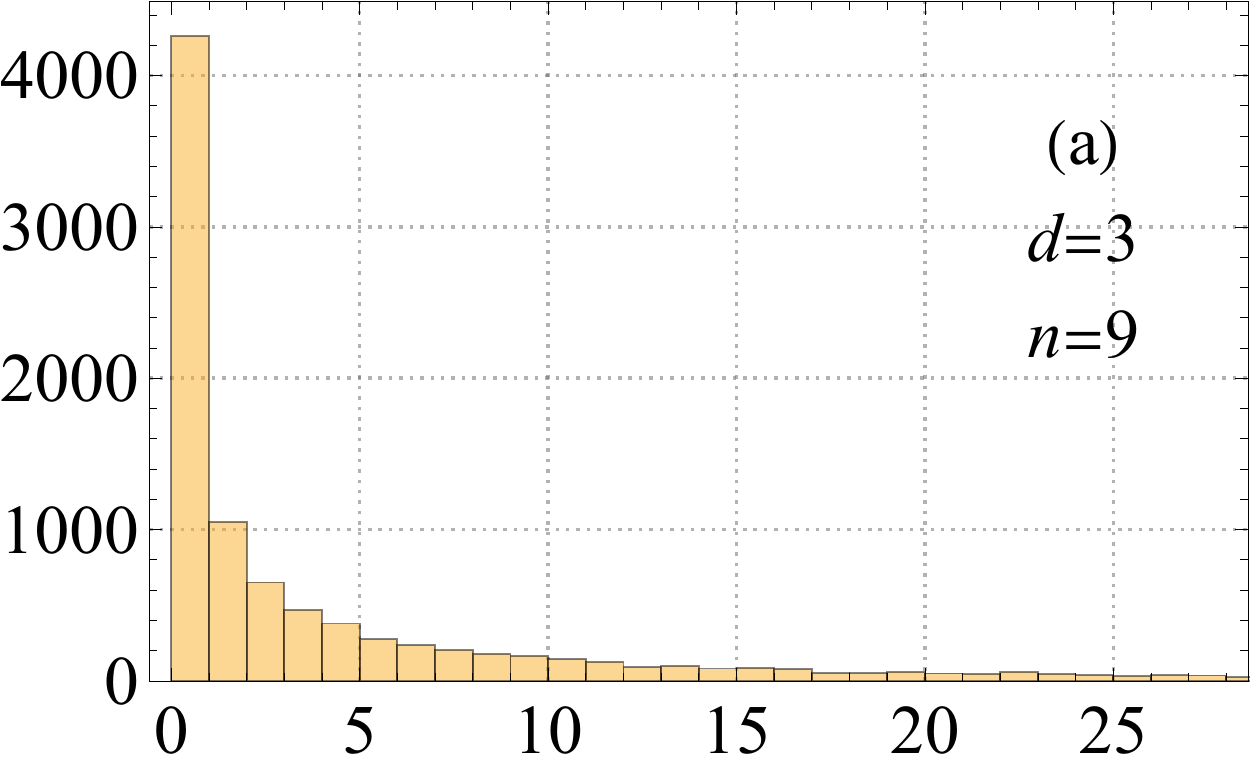}}}
\subfloat[$|\Lambda_{d,n}^{\rm best}-\Lambda_{d,n}^{\rm MaxMinSv}|$]{{\includegraphics[width=.33 \textwidth]{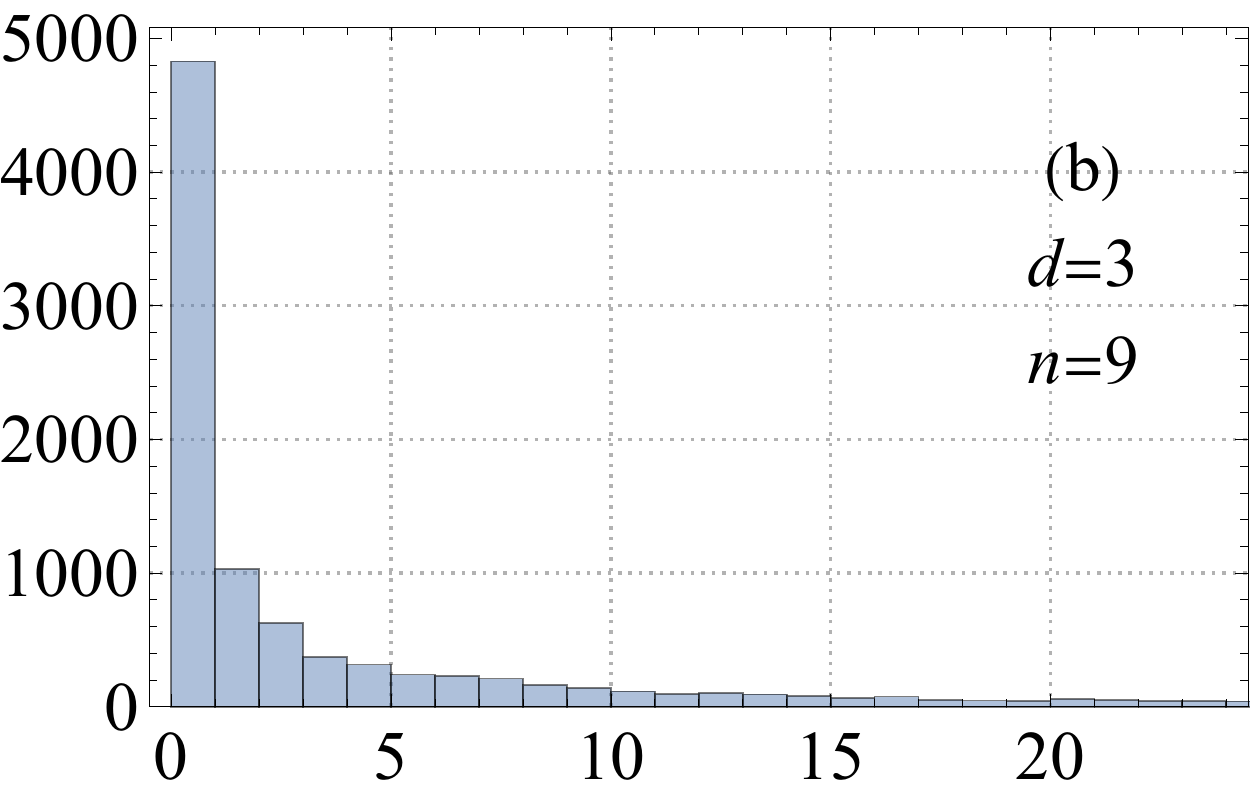}}}
\subfloat[$|\Lambda_{d,n}^{\rm MaxVol}-\Lambda_{d,n}^{\rm MaxMinSv}|$]{{\includegraphics[width=.33 \textwidth]{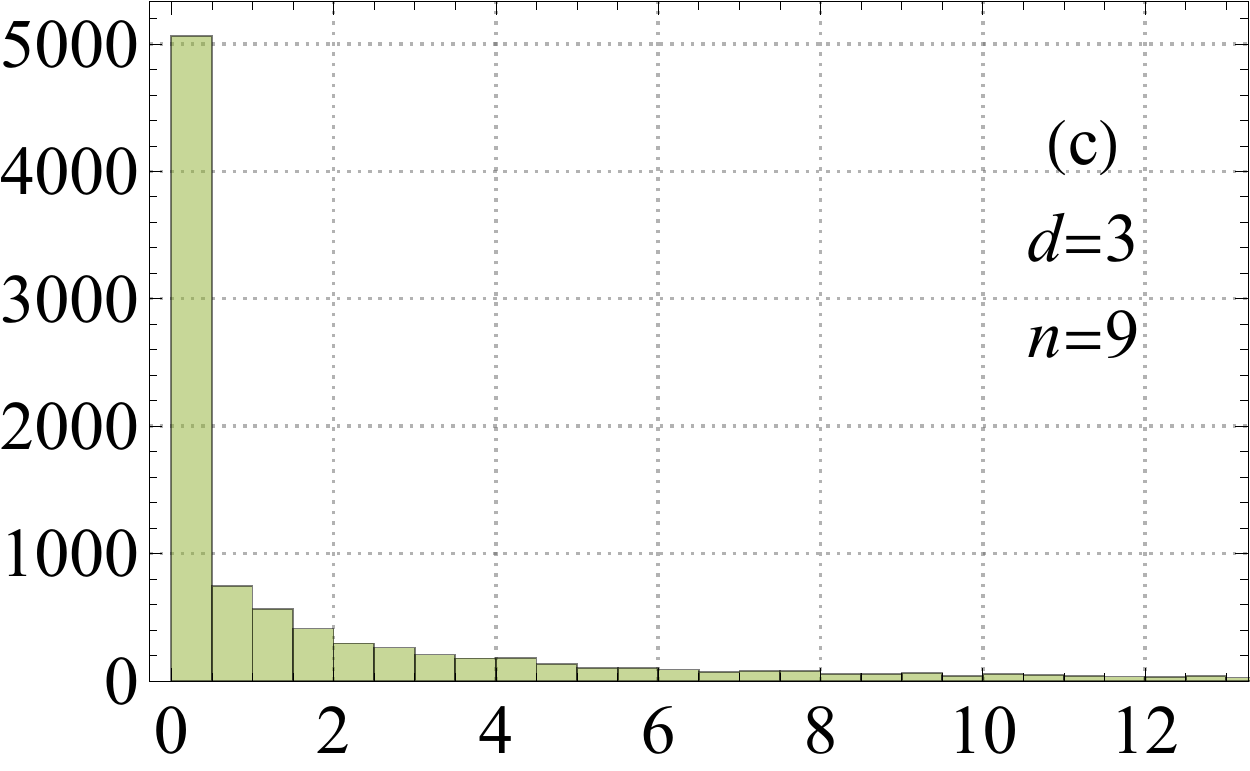}}}
\caption{Case (iii), where $10\,000$ realizations of uniformly random node sequences with $d=3$ were generated for each $n\in\{5,6,7,8,9\}$.  The occurrences of the differences (a)--(c) between the Lebesgue constants corresponding to differently chosen polynomial bases have been tabulated in the histograms above.}\label{testiii}
\end{figure}

In the first experiment, $10\,000$ different realizations of uniformly random nodes $\mathcal{X}=(\mathbf{x}_i)_{i=1}^n$ were generated in $[0,1]^d$ in three separate cases:
\begin{itemize}
\item[(i)] For $d=2$ and for each $n\in\{4,5\}$ with $\mathscr{B}=(\mathbf{x}^{\boldsymbol{\alpha}})_{|\boldsymbol{\alpha}|\leq 2}$;
\item[(ii)] For $d=2$ and for each $n\in\{7,8,9\}$ with $\mathscr{B}=(\mathbf{x}^{\boldsymbol{\alpha}})_{|\boldsymbol{\alpha}|\leq 3}$;
\item[(iii)] For $d=3$ and for each $n\in\{5,6,7,8,9\}$ with $\mathscr{B}=(\mathbf{x}^{\boldsymbol{\alpha}})_{|\boldsymbol{\alpha}|\leq 2}$.
\end{itemize}
For each realization of the node sequence $\mathcal{X}$, the generalized Vandermonde matrix $\mathscr{V}_{\mathscr{B},\mathcal{X}}$ can be constructed. Due to the low dimensionality of this example, all possible polynomial bases $\mathcal{B}\subset\mathscr{B}$ with $\#\mathcal{B}=n$ can be determined and out of these, three different polynomial bases were selected:
\begin{itemize}
\item The polynomial basis $\mathcal{B}_{\rm best}\subset\mathscr{B}$ producing the smallest Lebesgue constant denoted by $\Lambda_{d,n}^{\rm best}=\Lambda_{\mathcal{B}_{\rm best},\mathcal{X}}$.
\item The polynomial basis $\mathcal{B}_{\rm MaxVol}\subset\mathscr{B}$ corresponding to the maximum volume Vandermonde submatrix $V_{\mathcal{B}_{\rm MaxVol},\mathcal{X}}$ and its associated Lebesgue constant denoted by $\Lambda_{d,n}^{\rm MaxVol}=\Lambda_{\mathcal{B}_{\rm MaxVol},\mathcal{X}}$.
\item The polynomial basis $\mathcal{B}_{\rm MaxMinSv}\subset\mathscr{B}$ corresponding to the maximal minimum singular value Vandermonde submatrix $V_{\mathcal{B}_{\rm MaxMinSv},\mathcal{X}}$ and its associated Lebesgue constant denoted by $\Lambda_{d,n}^{\rm MaxMinSv}=\Lambda_{\mathcal{B}_{\rm MaxMinSv},\mathcal{X}}$.
\end{itemize}
We compute the differences
\begin{itemize}
\item[(a)] $|\Lambda_{d,n}^{\rm best}-\Lambda_{d,n}^{\rm MaxVol}|$,
\item[(b)] $|\Lambda_{d,n}^{\rm best}-\Lambda_{d,n}^{\rm MaxMinSv}|$,
\item[(c)] $|\Lambda_{d,n}^{\rm MaxVol}-\Lambda_{d,n}^{\rm MaxMinSv}|$,
\end{itemize} 
for every realization of the node set in each subcase of (i), (ii), and (iii), respectively. The occurrences of the differences (a)--(c) have been tabulated in the histograms displayed in Figures~\ref{testi}--\ref{testiii}.

In the case (iii) corresponding to $d=3$, several node configurations were encountered where either the maximum volume submatrix or the maximal minimum singular value submatrix were nearly singular, i.e., the determinant and minimum singular value were effectively zero; these cases have been dismissed from consideration. For $n$ equal to $5,6,7,8,$ and $9$, the numbers of dismissed realizations were $47$, $80$, $137$, $204$, and $241$, respectively. These cases can be handled by expanding the trial basis function set $\mathscr{B}$ with higher degree monomial basis functions or using basis functions other than monomials, but this was not done for this demonstration.

\begin{figure}[!h]
\captionsetup[subfigure]{labelformat=empty}
\centering
\subfloat{{\includegraphics[width=.33 \textwidth]{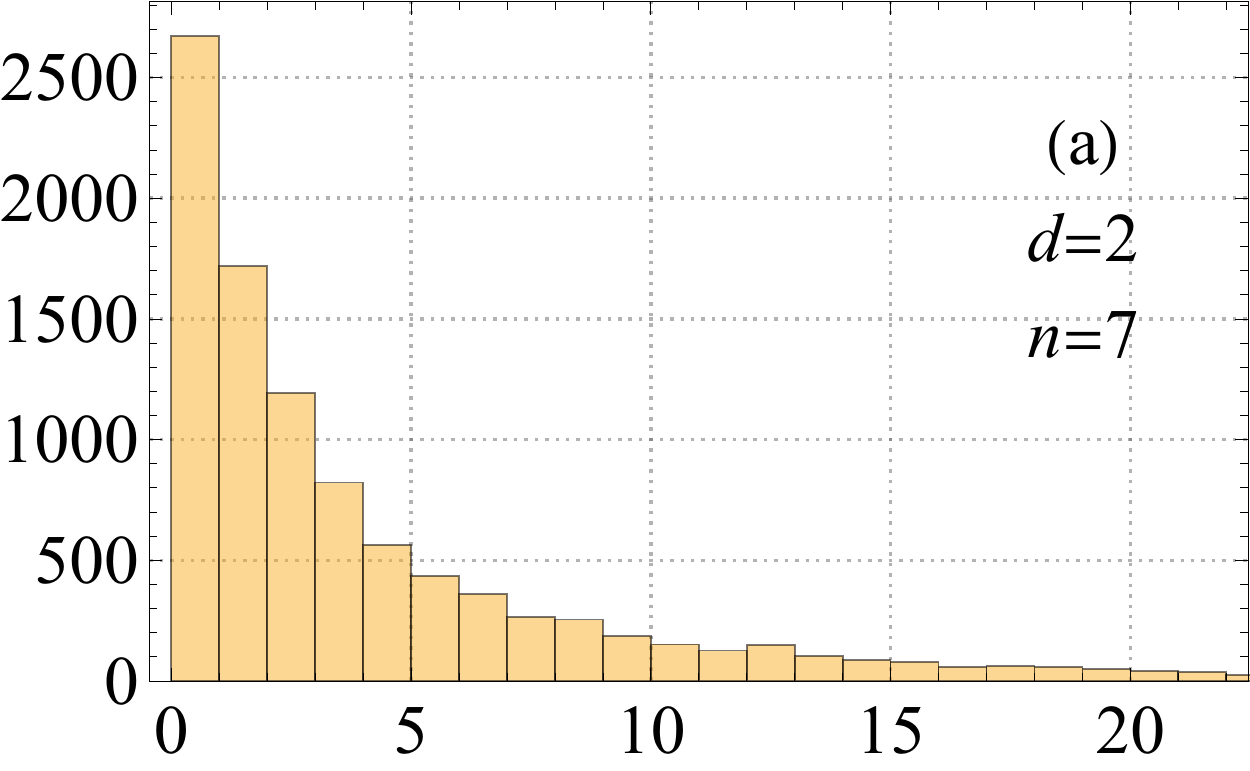}}}
\subfloat{{\includegraphics[width=.33 \textwidth]{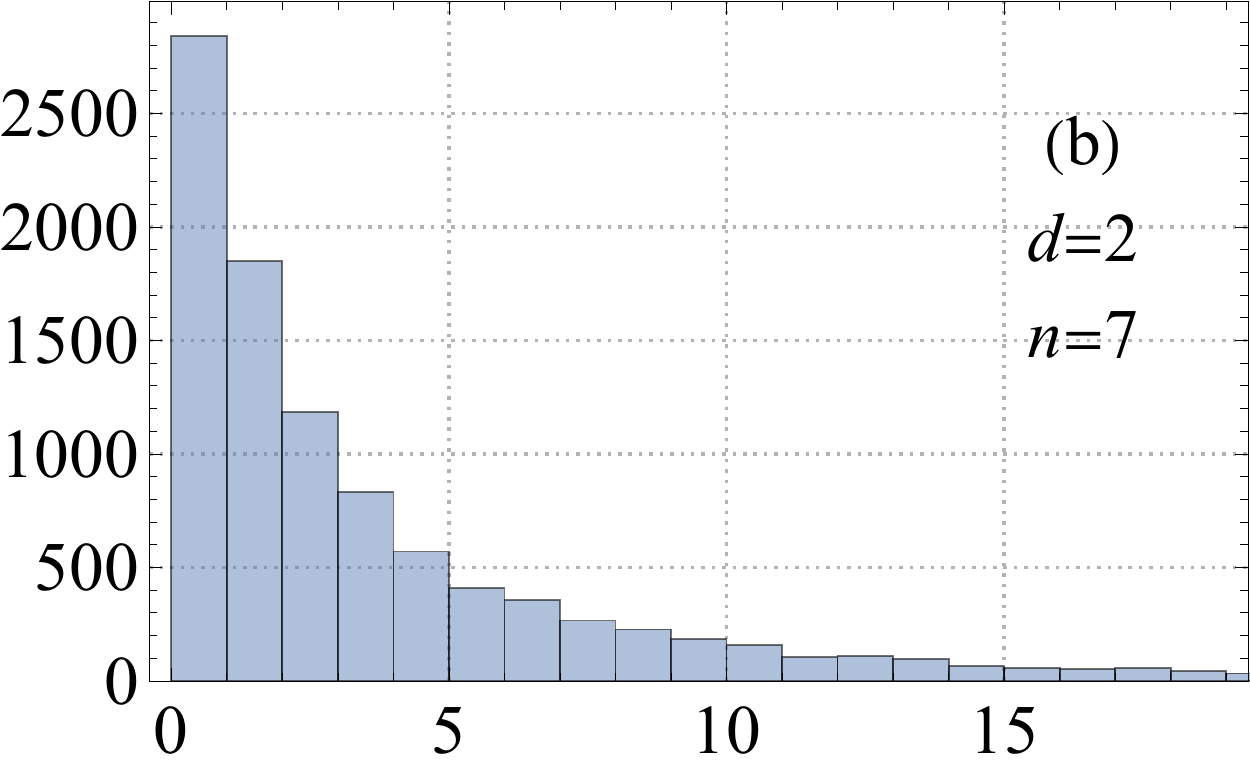}}}
\subfloat{{\includegraphics[width=.33 \textwidth]{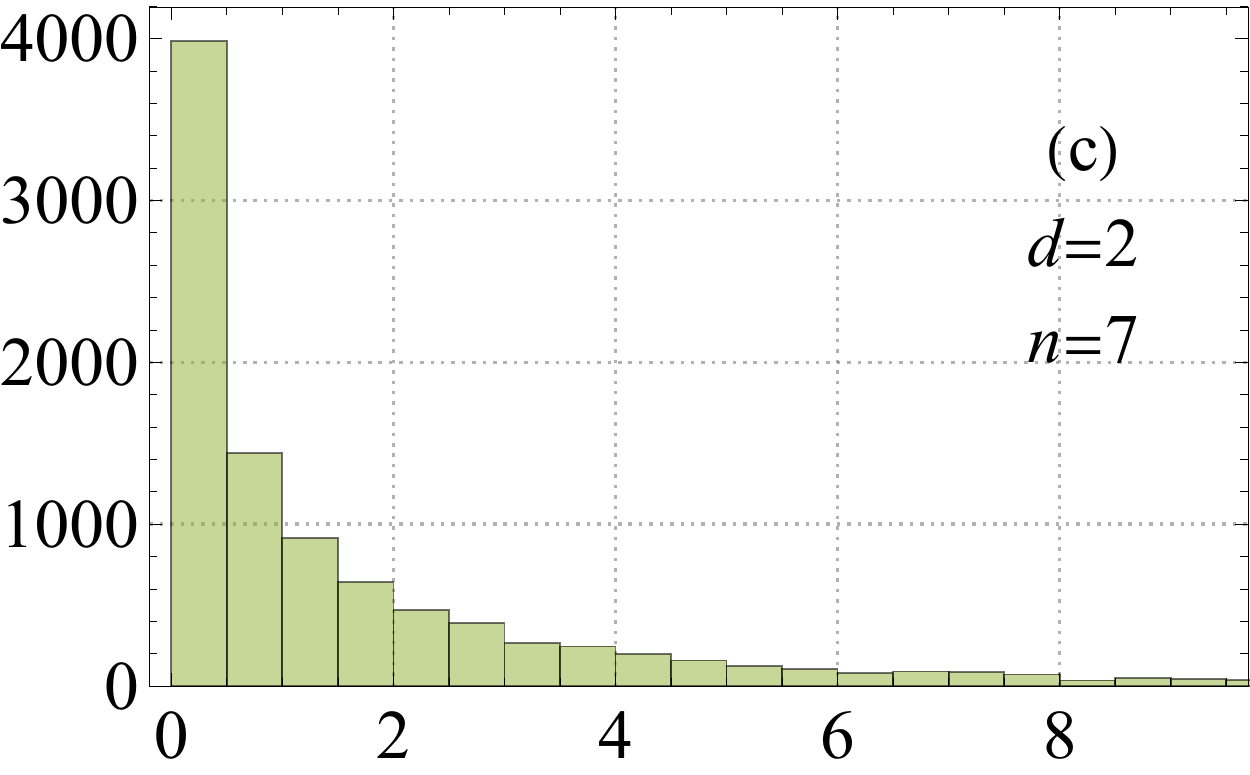}}}\\
\subfloat{{\includegraphics[width=.33 \textwidth]{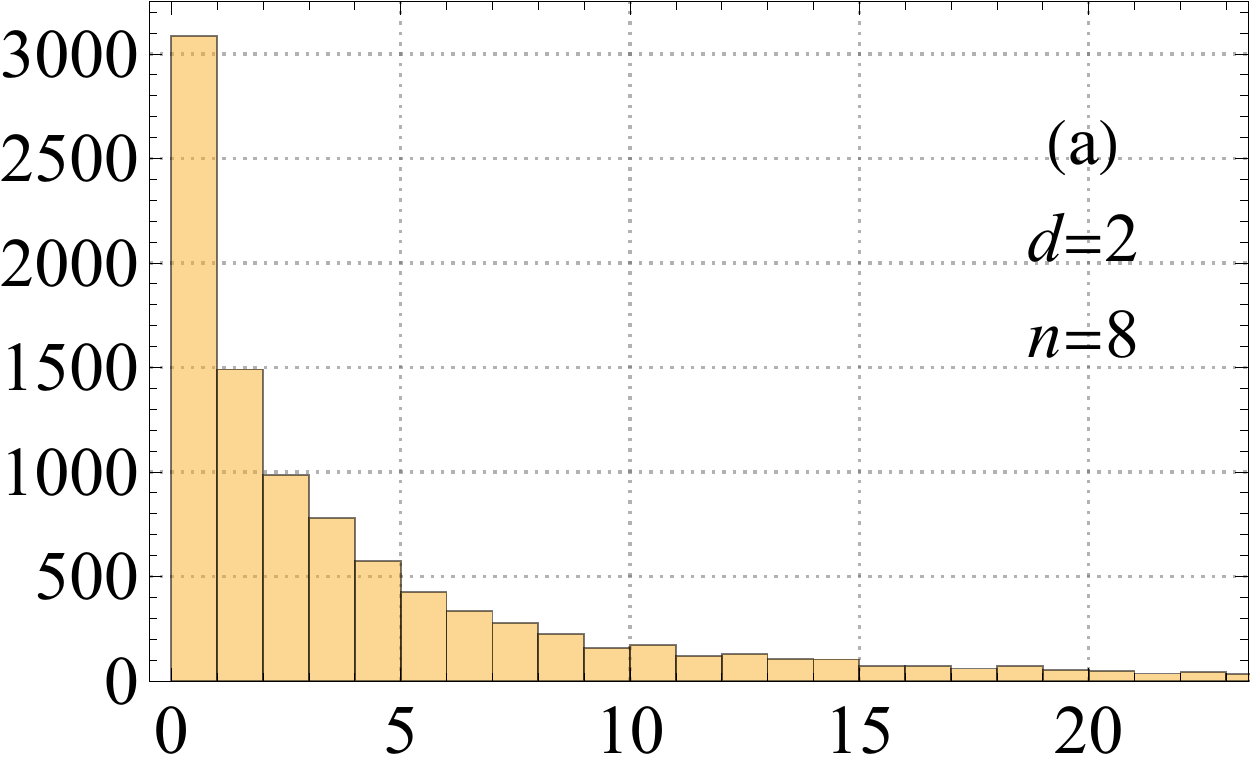}}}
\subfloat{{\includegraphics[width=.33 \textwidth]{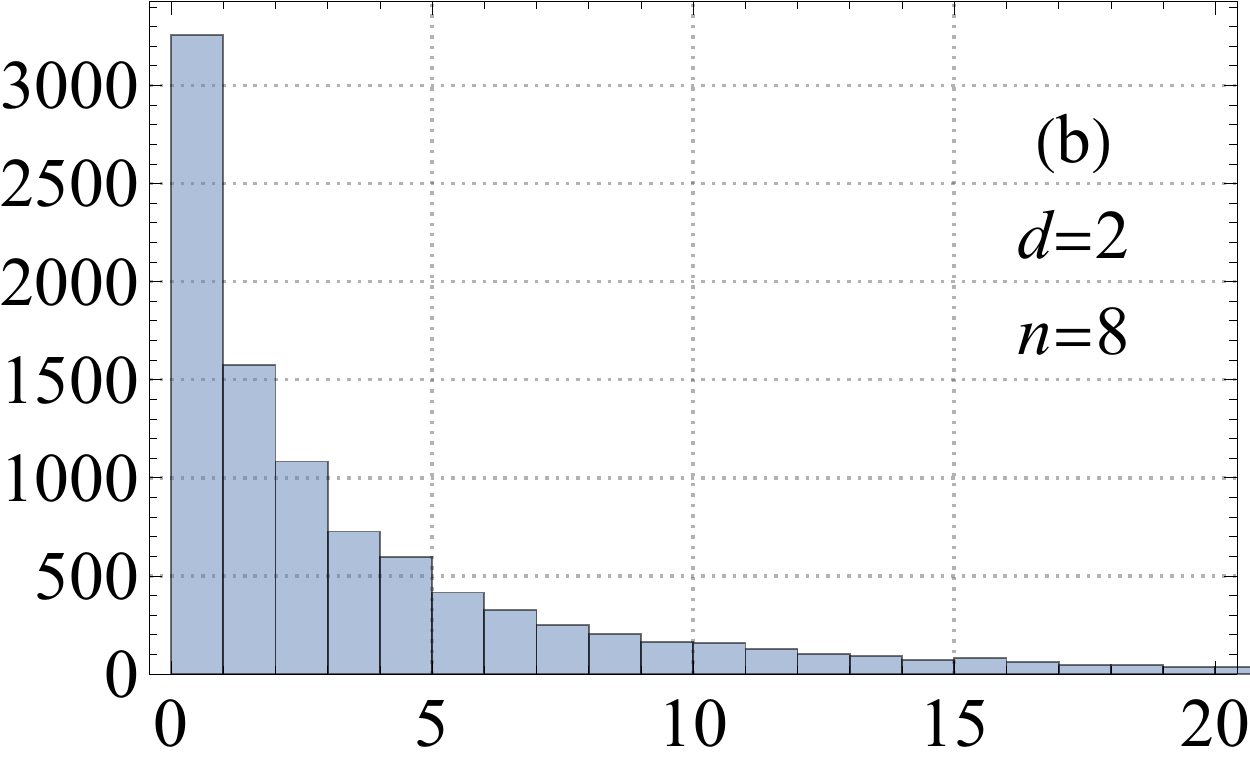}}}
\subfloat{{\includegraphics[width=.33 \textwidth]{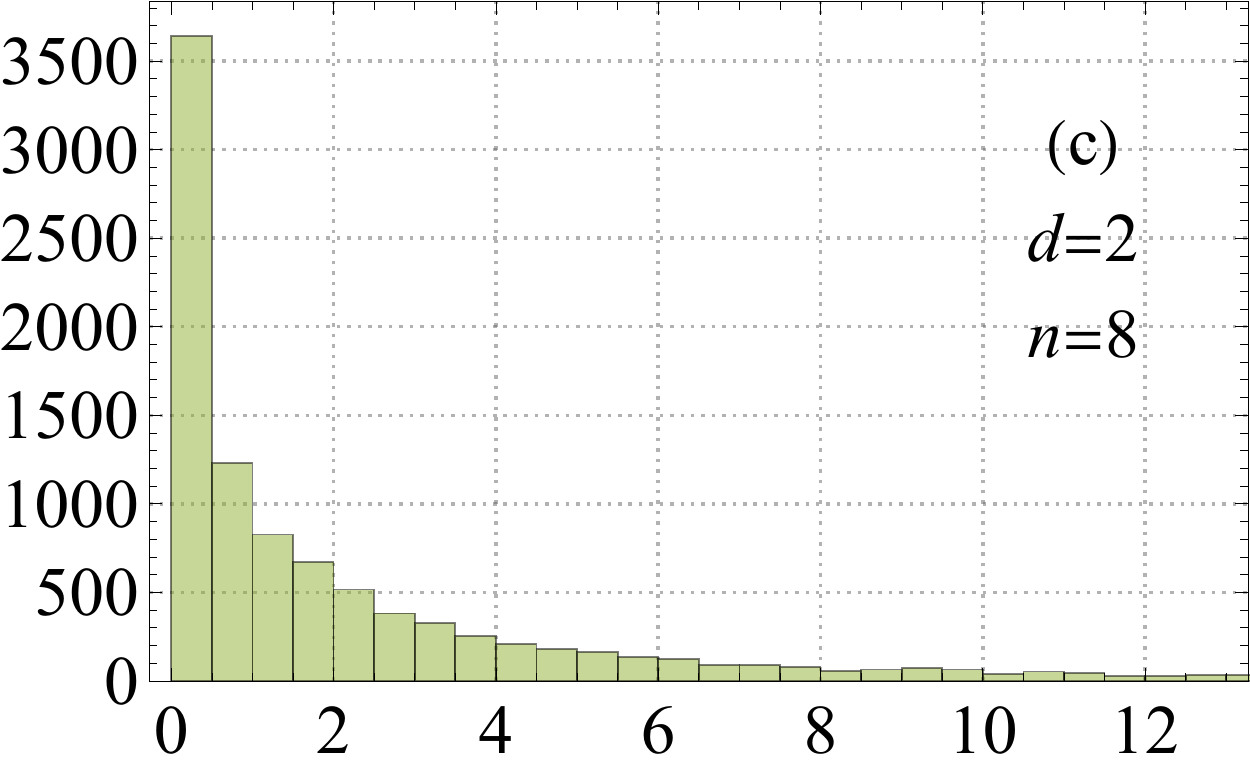}}}\\
\subfloat[$|\Lambda_{n,d}^{\rm best}-\Lambda_{n,d}^{\rm MaxVol}|$]{{\includegraphics[width=.33 \textwidth]{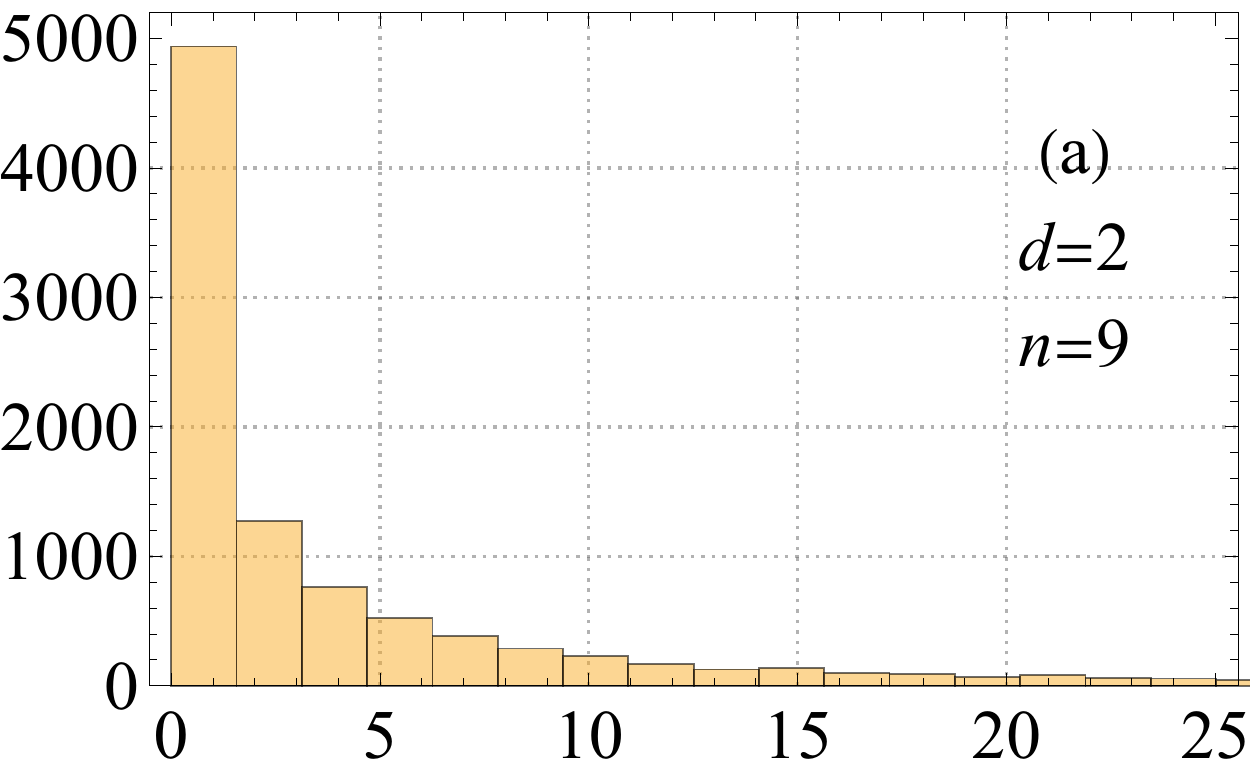}}}
\subfloat[$|\Lambda_{n,d}^{\rm best}-\Lambda_{n,d}^{\rm MaxMinSv}|$]{{\includegraphics[width=.33 \textwidth]{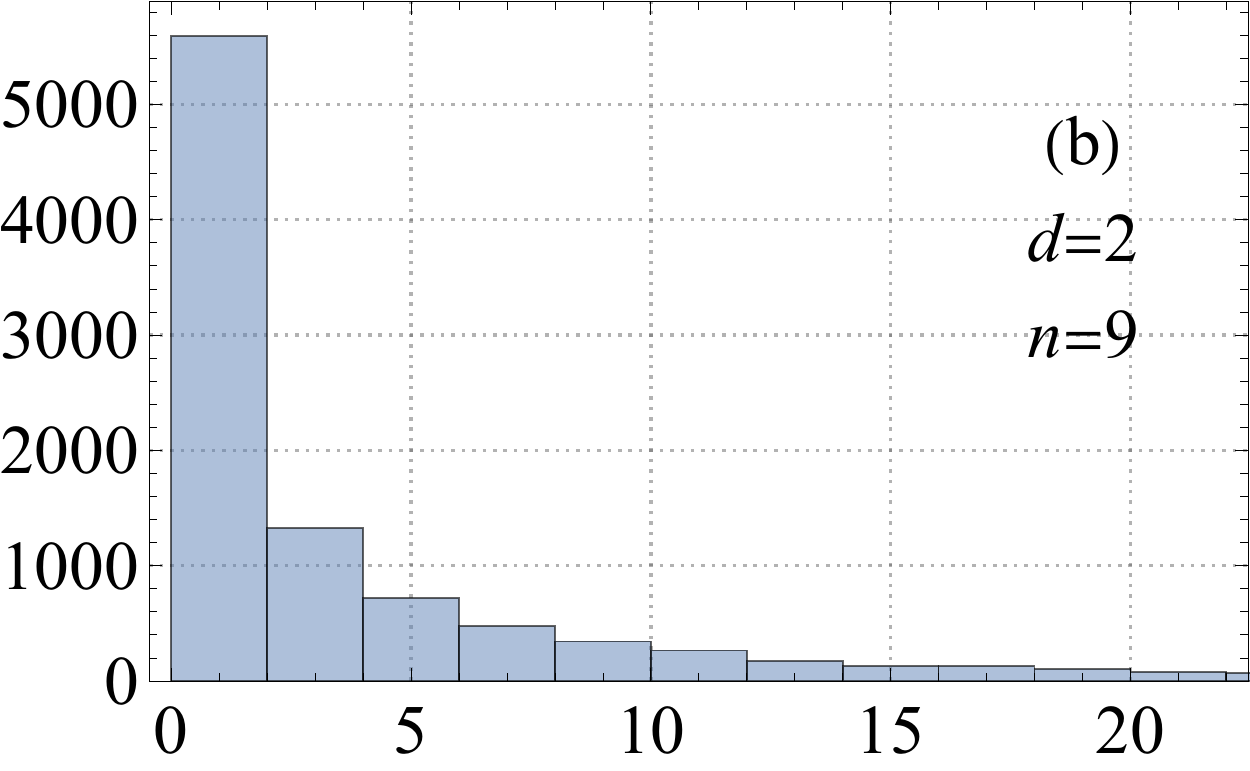}}}
\subfloat[$|\Lambda_{n,d}^{\rm MaxVol}-\Lambda_{n,d}^{\rm MaxMinSv}|$]{{\includegraphics[width=.33 \textwidth]{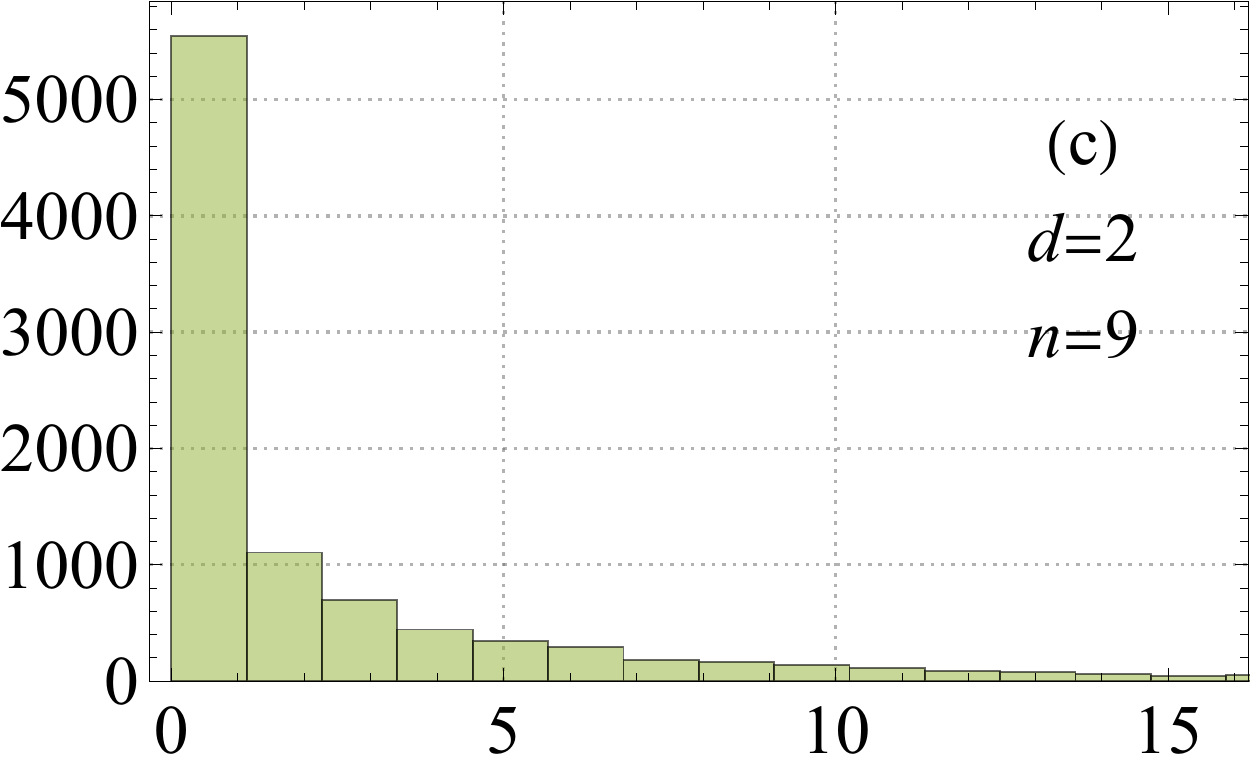}}}
\caption{Case (ii), where $10\,000$ realizations of uniformly random node sequences with $d=2$ were generated for each $n\in\{7,8,9\}$.  The occurrences of the differences (a)--(c) between the Lebesgue constants corresponding to differently chosen polynomial bases have been tabulated in the histograms above.}
\end{figure}

From the numerical experiments on uniformly random nodes, we find that neither the maximum volume submatrix nor the maximal minimum singular value submatrix guarantee obtaining the basis with optimal Lebesgue constant for every trial. However, the obtained bases produce in most cases a Lebesgue constant that is very close to optimum if not optimal, and the occurrences with a large difference to optimum are extremely rare---hence the exponentially vanishing tail in the histograms. When $n$ approaches the dimension of a total degree polynomial space, the frequency of obtaining a near-optimal Lebesgue constant with either approach increases. Moreover, there do not appear to be notable differences in the performance of the maximum volume and maximal minimum singular value Vandermonde submatrices. This is a welcome observation since the former type is easier to approximate due to the availability of applicable numerical algorithms.

\subsection{Incomplete sparse grids}\label{snodes}

\subsubsection{Smolyak interpolating polynomial}

The Smolyak interpolating polynomial is a well-known tool used to extend univariate interpolation rules defined on a given interval $\mathcal{I}\subset\mathbb{R}$ to the hypercube $\mathcal{I}^d\subset\mathbb{R}^d$ in a computationally efficient manner. We refer to~\cite{barthelmann00} for a detailed account on the construction of the Smolyak polynomial and simply give the formula for the order $k\geq 0$ Smolyak polynomial in $d$ variables in the special case where it admits to a general expression in the form
\[
\mathcal{S}_{d,k}f(\mathbf{x})=\sum_{\substack{\boldsymbol{\alpha}\in\mathbb{Z}_+^d\\ d\leq |\boldsymbol{\alpha}|\leq d+k}}\sum_{i_1=m(\alpha_1-1)+1}^{m(\alpha_1)}\cdots \sum_{i_d=m(\alpha_d-1)+1}^{m(\alpha_d)}c_{i_1,\ldots,i_d}(f)T_{i_1-1}(x_1)\cdots T_{i_d-1}(x_d),
\]
where $\mathbf{x}=[x_1,\ldots,x_d]^\textup{T}\in\mathbb{R}^d$, we set $m(0)=0$, $m(1)=1$, and $m(k)=2^{k-1}+1$ for $k>1$, respectively, and the coefficients $c_{i_1,\ldots,i_d}(f)\in\mathbb{R}$ are determined uniquely by the data we wish to interpolate. Here, the univariate polynomials $T_k$ are Chebyshev polynomials of the first kind defined by the three-term recursion
\[
T_0(x)=1,\quad T_1(x)=x,\quad \text{and}\quad T_{k+1}(x)=2xT_k(x)-T_{k-1}(x)\quad \text{for }k\geq 1.
\]
Let us denote the univariate basis $B_{k+1}=(T_j)_{j=0}^k$ for $k\geq 0$ and let us also denote the sets of univariate Clenshaw--Curtis abscissae in the interval $[-1,1]$ by setting
\[
X_1^{\rm CC}=\{0\}\quad\text{and}\quad X_k^{\rm CC}=\left\{-\cos\left(\frac{(j-1)\pi}{m(k)-1}\right):j=1,\ldots,m(k)\right\}\quad\text{for } k>1.
\]
The Smolyak interpolating polynomial is unisolvent with respect to the basis
\[
\mathcal{B}_{d,k}=\bigcup_{\substack{\boldsymbol{\alpha}\in\mathbb{Z}_+^d\\ d\leq|\boldsymbol{\alpha}|\leq d+k}}B_{m(\alpha_1)}\times \cdots\times B_{m(\alpha_d)},\quad k\geq 0,
\]
and the interpolating polynomial is exact on the nodes that form the \emph{sparse grid}
\[
\mathcal{X}_{d,k}=\bigcup_{\substack{\boldsymbol{\alpha}\in\mathbb{Z}_+^d\\ |\boldsymbol{\alpha}|=d+k}}X_{\alpha_1}^{\rm CC}\times\cdots\times X_{\alpha_d}^{\rm CC},\quad k\geq 0.
\]

The Smolyak interpolating polynomial $\mathcal{S}_{d,k}f$ can be obtained for the respective basis $\mathcal{B}_{d,k}$, node configuration $\mathcal{X}_{d,k}=\{\mathbf{x}_1,\ldots,\mathbf{x}_n\}$, and data $\mathbf{f}=[f(\mathbf{x}_1),\ldots,f(\mathbf{x}_n)]^\textup{T}$ by solving the coefficients $c_{i_1,\ldots,i_d}(f)$ implicitly from the Vandermonde system
\[
V_{\mathcal{B}_k,\mathcal{X}_k}^\textup{T} \mathbf{c}=\mathbf{f}
\]
which yields the Smolyak interpolating polynomial
\[
\mathcal{S}_{d,k}f(\mathbf{x})=\sum_{i=1}^n c_i\phi_i(\mathbf{x}),\quad \phi_i\in\mathcal{B}_{d,k},
\]
and satisfies
\[
\mathcal{S}_{d,k}f(\mathbf{x}_i)=f(\mathbf{x}_i)\quad\text{for all }\mathbf{x}_i\in\mathcal{X}_{d,k}.
\]

The Smolyak interpolating polynomial is known to be very stable. However, it is only well-defined over \emph{complete} sparse grids $\mathcal{X}_{d,k}$. In the following, we investigate the stability of the interpolating polynomial when one tries to interpolate over \emph{incomplete} sparse grid node sets $\mathcal{Y}$ that lie between two complete sparse grids such that $\mathcal{X}_{d,k}\subsetneq\mathcal{Y}\subsetneq\mathcal{X}_{d,k+1}$. We use the \texttt{MaxVol} algorithm~\cite{maxvol} on the rows of the generalized Vandermonde matrix $\mathscr{V}_{\mathcal{B}_{d,k+1},\mathcal{Y}}$ to find a basis $\mathcal{B}_{\rm MaxVol}\subset\mathcal{B}_{d,k+1}$ for which the Lagrange interpolation problem~\eqref{lagrangeproblem} is well-defined over $\mathcal{Y}$.

\subsubsection{Incomplete Smolyak interpolating polynomial}
We denote $n_k=\#\mathcal{X}_{d,k}$ and suppress the dependence on $d$  since there is no risk of confusion in the sequel. We consider in the following an incrementally increasing sequence of nodes $\mathcal{Y}_{d,k,i}$ with $\#\mathcal{Y}_{d,k,i}=\#\mathcal{X}_{d,k}+i$ containing nodes which lie between successive sparse grids $\mathcal{X}_{d,k}$ and $\mathcal{X}_{d,k+1}$ such that
\[
\mathcal{X}_{d,k}\subsetneq \mathcal{Y}_{d,k,1}\subsetneq\cdots\subsetneq\mathcal{Y}_{d,k,n_{k+1}-n_k-1}\subsetneq\mathcal{X}_{d,k+1},
\]
and measure the stability and accuracy of the resulting interpolating polynomials over incomplete sparse grids by computing their associated Lebesgue constants.

\begin{figure}[!h]
\centering
\subfloat[The ordering of the complete sparse grid node configuration $\mathcal{X}_{2,2}=\{\mathbf{x}_1,\ldots,\mathbf{x}_{13}\}$.]{{\includegraphics[width=.31 \textwidth]{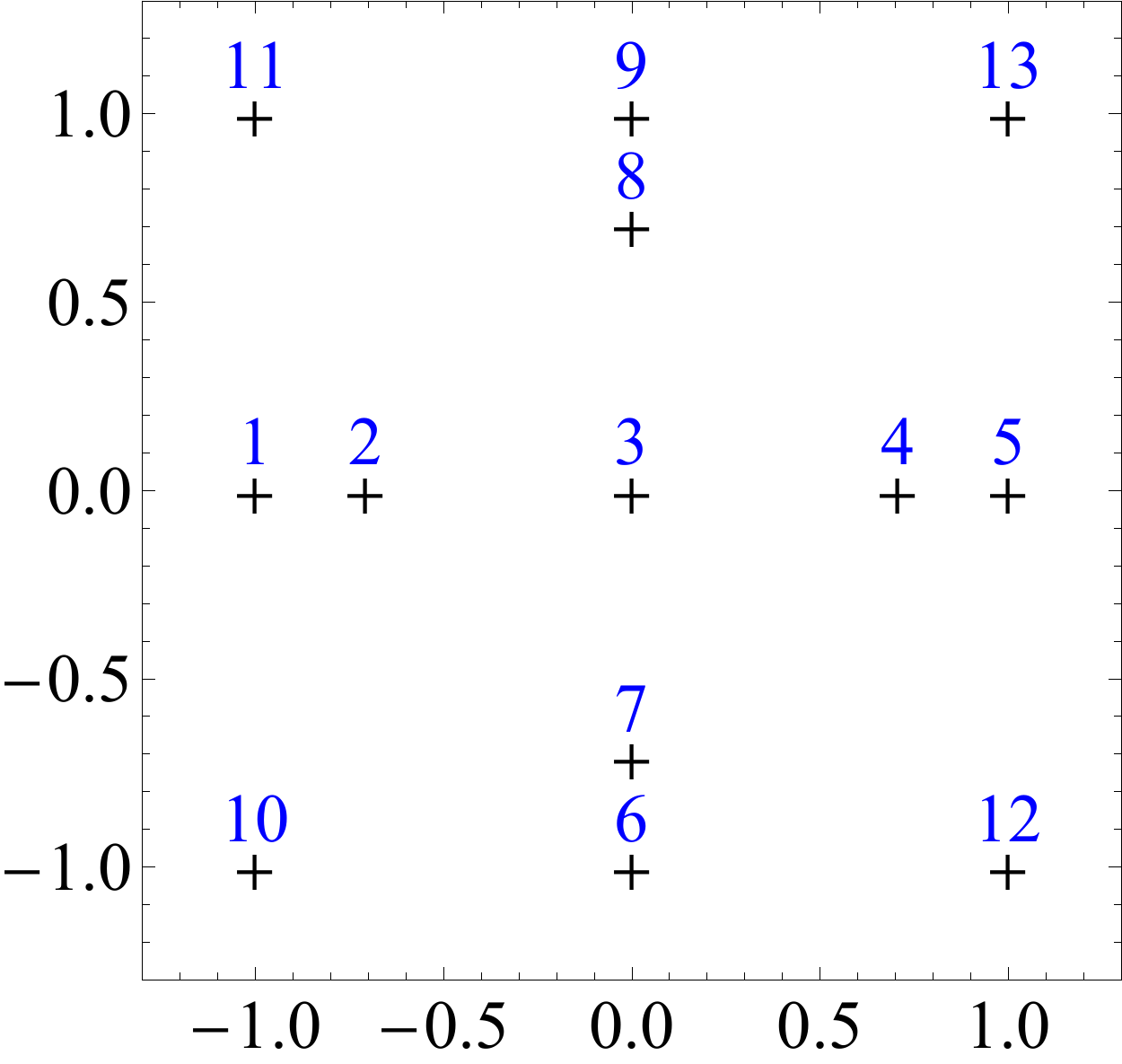}}}\quad 
\subfloat[The ordering of the incomplete sparse grid node configuration $\mathcal{X}_{2,3}\setminus\mathcal{X}_{2,2}=\{\mathbf{x}_{14},\ldots,\mathbf{x}_{29}\}.$]{{\includegraphics[width=.31 \textwidth]{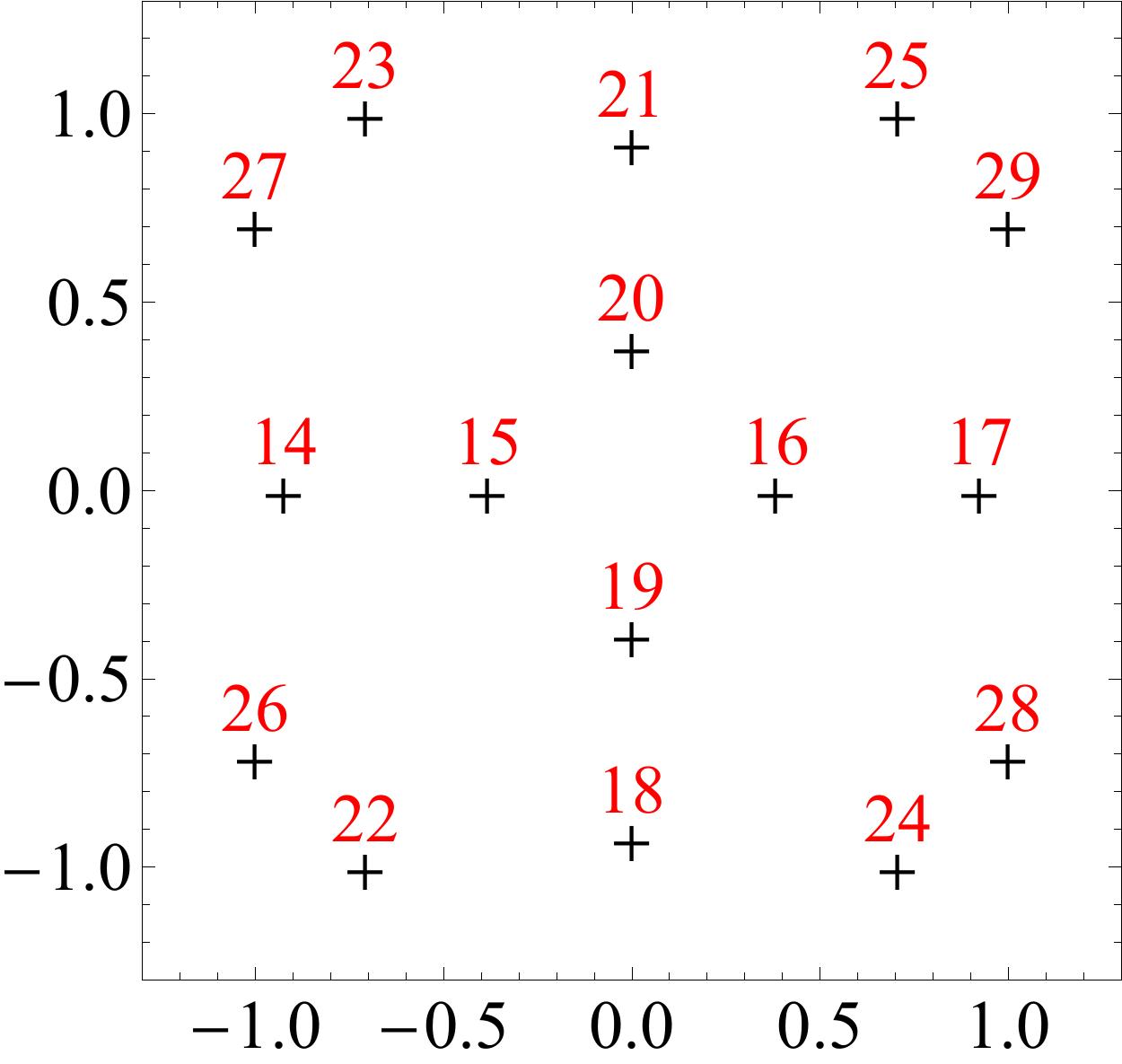}}}\quad 
\subfloat[The complete sparse grid node configuration $\mathcal{X}_{2,3}=\{\mathbf{x}_1,\ldots,\mathbf{x}_{29}\}$.]{{\includegraphics[width=.31 \textwidth]{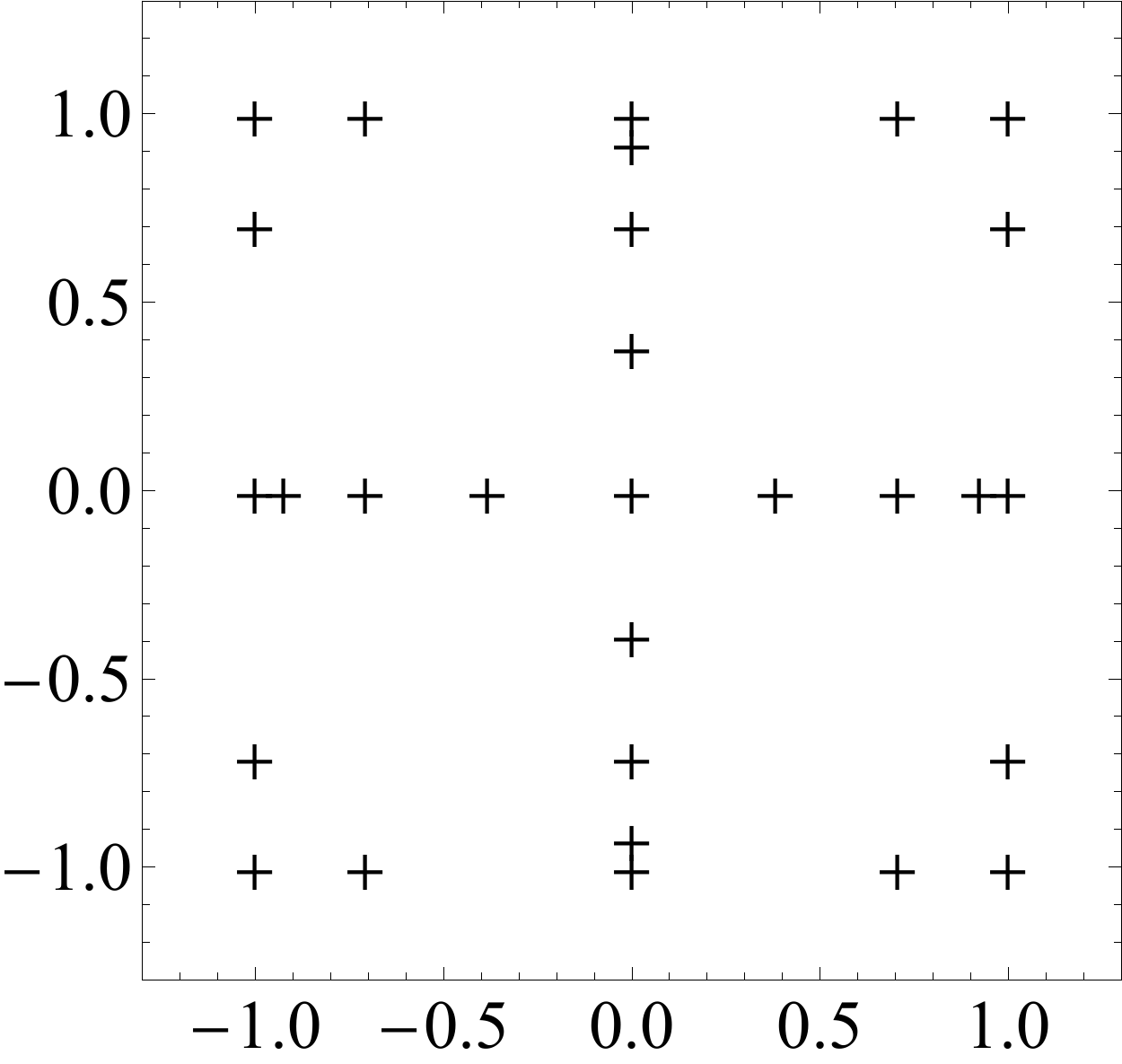}}}
\caption{The node ordering used for the Smolyak--Clenshaw--Curtis interpolation rule in the case (iv) for $d=2$.}\label{orderingfig}
\end{figure}

Let us denote the elements of the sparse grids $\mathcal{X}_{d,k}$ and incomplete sparse grids $\mathcal{Y}_{d,k,i}$ by
\begin{align*}
&\mathcal{X}_{d,k}=\{\mathbf{x}_1,\ldots,\mathbf{x}_{n_k}\},\\
&\mathcal{Y}_{d,k,i}=\{\mathbf{x}_1,\ldots,\mathbf{x}_{n_k+i}\},\quad i\in\{1,\ldots,n_{k+1}-n_k-1\}.
\end{align*}
We consider the following cases:
\begin{itemize}
\item[(iv)] The Lebesgue constants of the incomplete sparse grids $\mathcal{X}_{d,2}\subsetneq \mathcal{Y}_{d,2,1}\subsetneq\cdots\subsetneq\mathcal{Y}_{d,2,n_3-n_2+1}\subsetneq\mathcal{X}_{d,3}$ for $d\in\{2,3\}$.
\item[(v)] The Lebesgue constants of the incomplete sparse grids $\mathcal{X}_{d,3}\subsetneq \mathcal{Y}_{d,3,1}\subsetneq\cdots\subsetneq\mathcal{Y}_{d,3,n_4-n_3+1}\subsetneq\mathcal{X}_{d,4}$ for $d\in\{2,3\}$.
\end{itemize}
For the purpose of demonstration, we fix the ordering of the nodes. The ordering in the case (iv) for $d=2$ is illustrated in Figure~\ref{orderingfig}; the tables containing the explicit enumeration of the nodes and basis functions for both subcases of (iv) are given explicitly in Appendix~\ref{sappendix}. In the case (v), the node ordering is produced by using the same algorithm for sparse grid generation as in the case (iv), but the explicit numbering is omitted for brevity.

\begin{figure}[!h]
\captionsetup[subfigure]{labelformat=empty}
\centering
\subfloat[The case (iv) for $d=2$.]{{\includegraphics[width=.495 \textwidth]{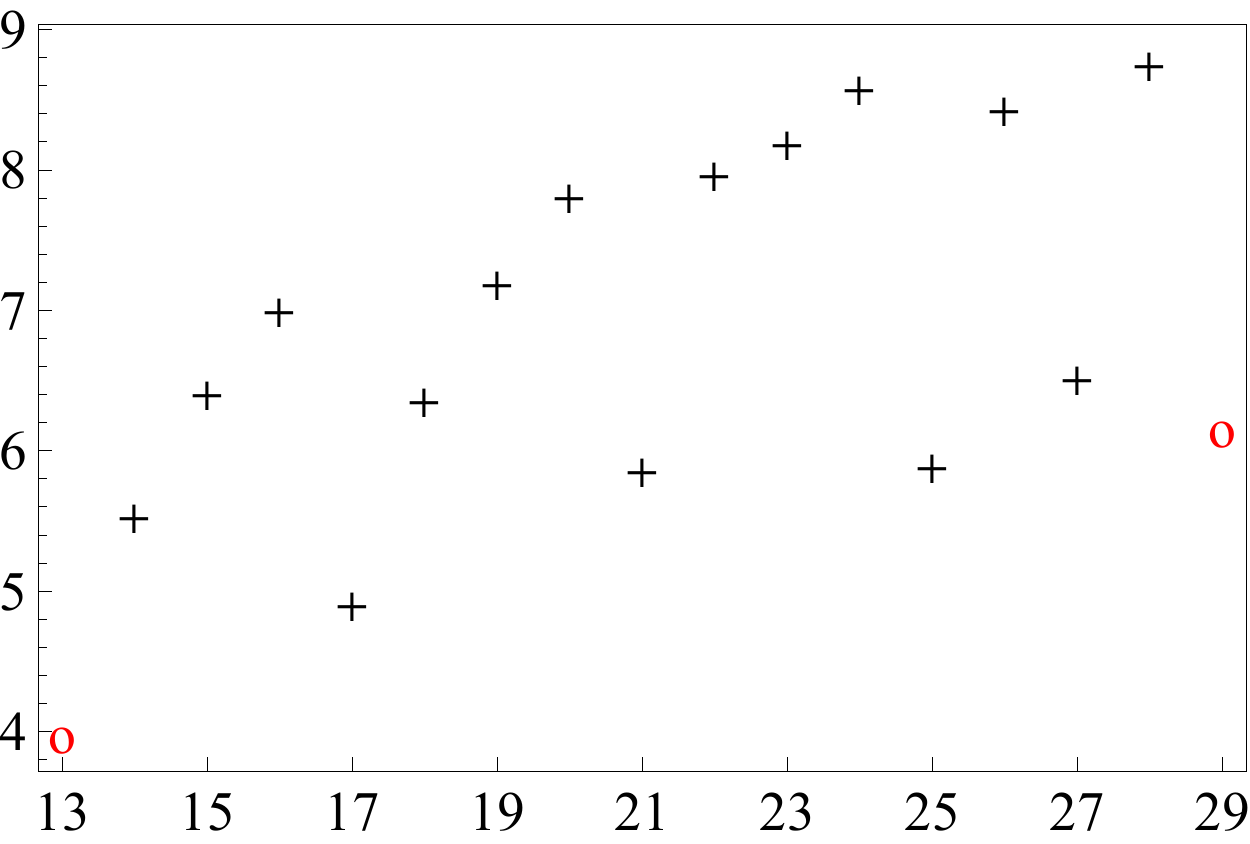}}}
\subfloat[The case (iv) for $d=3$.]{{\includegraphics[width=.501 \textwidth]{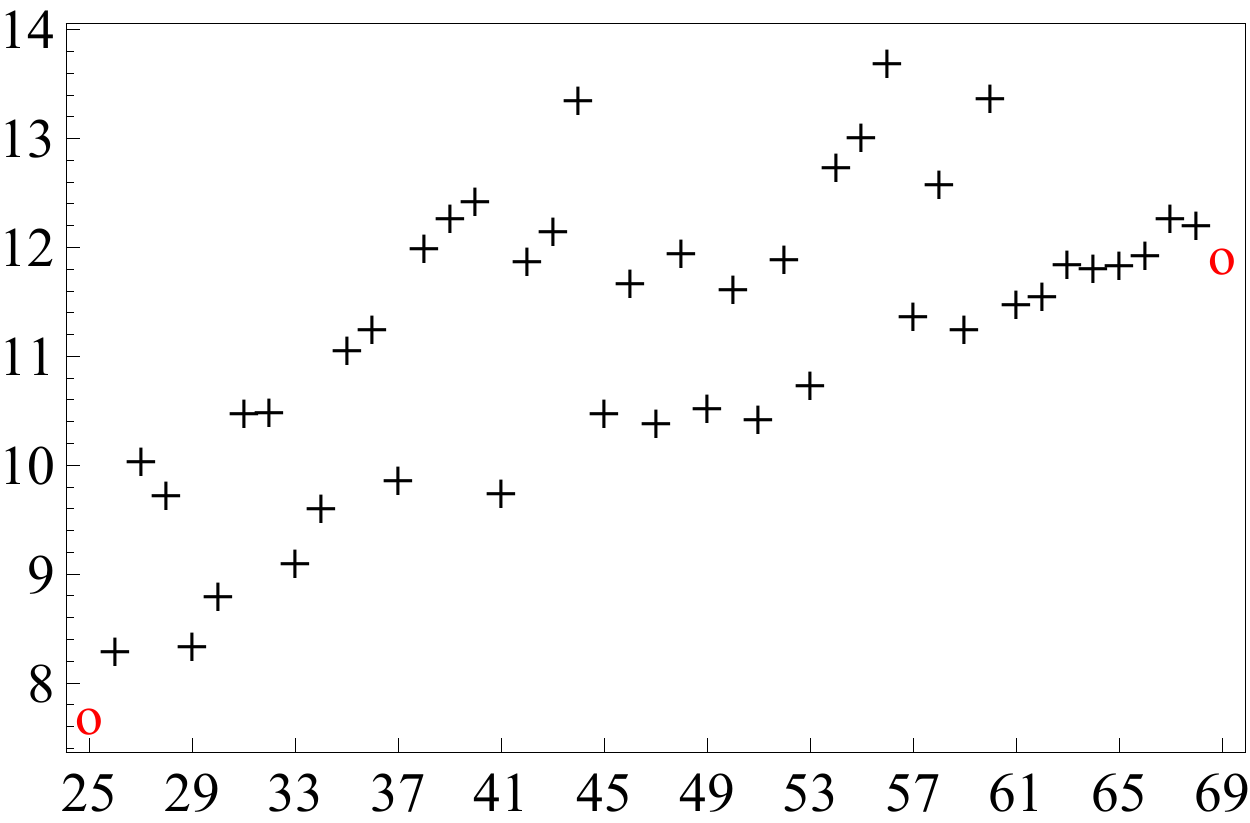}}}
\caption{The obtained Lebesgue constants for polynomial interpolation over incomplete sparse grids in case (iv). The horizontal axis corresponds to the cardinality of the incomplete sparse grid $\mathcal{Y}_{d,k,i}$. The circles represent the Lebesgue constants of the complete Smolyak interpolating polynomials for $k=2$ and $k=3$, respectively.}\label{sresults1}
\end{figure}

\begin{figure}[!h]
\captionsetup[subfigure]{labelformat=empty}
\centering
\subfloat[The case (v) for $d=2$.]{{\includegraphics[width=.495 \textwidth]{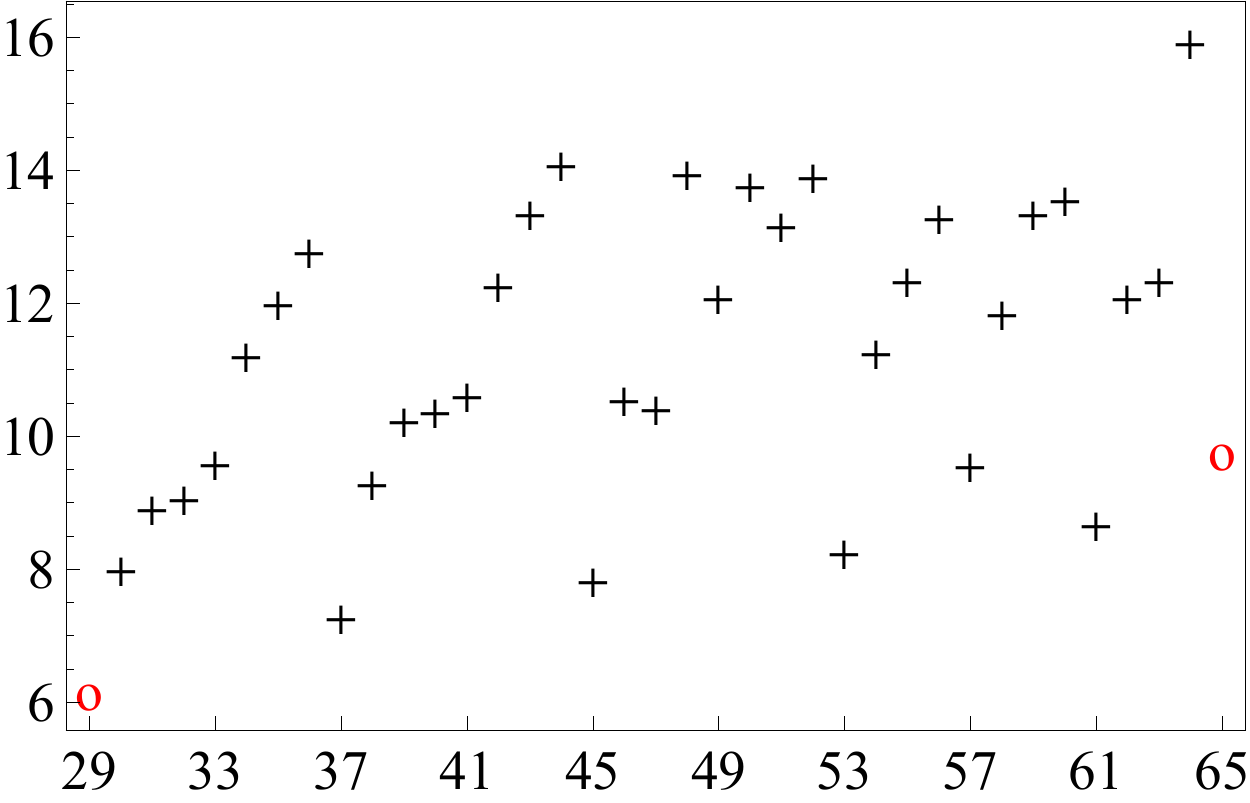}}}
\subfloat[The case (v) for $d=3$.]{{\includegraphics[width=.501 \textwidth]{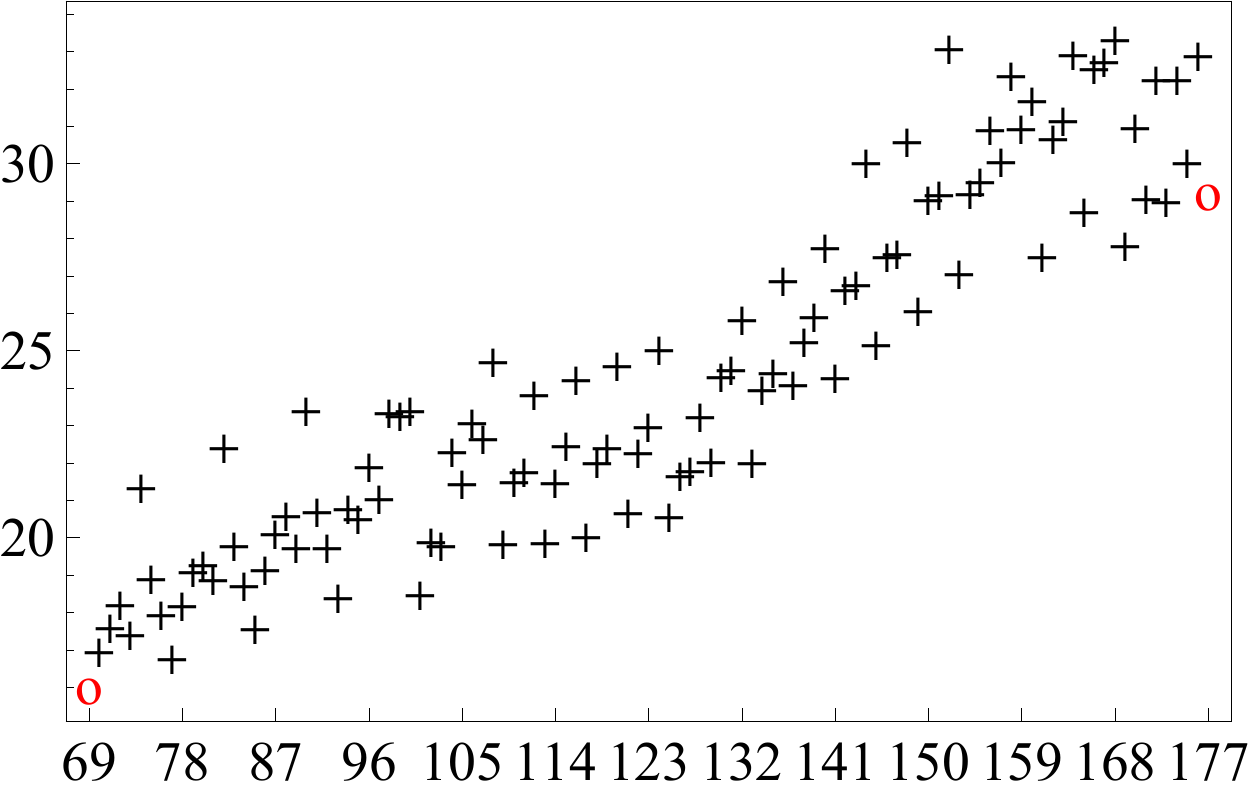}}}
\caption{The obtained Lebesgue constants for polynomial interpolation over incomplete sparse grids in case (v). The horizontal axis corresponds to the cardinality of the incomplete sparse grid $\mathcal{Y}_{d,k,i}$. The circles represent the Lebesgue constants of the complete Smolyak interpolating polynomials for $k=3$ and $k=4$, respectively.}\label{sresults2}
\end{figure}

To find the interpolating polynomial for each set $\mathcal{Y}_{d,k,i}$, we formulate the generalized Vandermonde matrix of the form
\[
\mathscr{V}_{\mathcal{B}_{d,k+1},\mathcal{Y}_{d,k,i}}=\begin{bmatrix}V_{\mathcal{B}_{d,k},\mathcal{X}_{d,k}}&\mathscr{V}_{\mathcal{B}_{d,k},\mathcal{Y}_{d,k,i}\setminus\mathcal{X}_{d,k}}\\ \mathscr{V}_{\mathcal{B}_{d,k+1}\setminus\mathcal{B}_{d,k},\mathcal{X}_{d,k}}&\mathscr{V}_{\mathcal{B}_{d,k+1}\setminus\mathcal{B}_{d,k},\mathcal{Y}_{d,k,i}\setminus\mathcal{X}_{d,k}}\end{bmatrix}.
\]
Then by using the \texttt{MaxVol} algorithm~\cite{maxvol} with respect to the rows of $\mathscr{V}_{\mathcal{B}_{d,k+1},\mathcal{Y}_{d,k,i}}$ we can find an approximate maximum volume submatrix
\[
V_{\mathcal{B}_{\rm MaxVol},\mathcal{Y}_{d,k,i}}
\]
the rows of which determine a basis $\mathcal{B}_{\rm MaxVol}\subset \mathcal{B}_{d,k+1}$ for which the Lagrange interpolation problem is well-defined with respect to $\mathcal{Y}_{d,k,i}$. We compute the Lebesgue constant over the convex hull $[-1,1]^d$ of $\mathcal{Y}_{d,k,i}$. The Lebesgue constants obtained for case (iv) are displayed in Figure~\ref{sresults1} and the results for case (v) are given in Figure~\ref{sresults2}.

Note that the ordering of the basis functions does not matter in this experiment. On the other hand, the ordering used for the node sequences in this example is a direct result of the algorithm used to generate the sparse grids: Changing the ordering of the nodes (symmetry notwithstanding) may change the obtained Lebesgue constants of the \emph{incomplete} sparse grids, but the results should be comparable as long as the numbering of the nodes is consistent between each multi-index set used in Smolyak's construction.

Notably, the Lebesgue constants of the incomplete sparse grids have magnitudes similar to those of the complete Smolyak interpolating polynomials. Moreover, the Lebesgue constants are at their smallest whenever a complete filament of the sparse grid is completed. For example, in the case (iv) for $d=2$ this corresponds to the grids $\mathcal{Y}_{2,2,i}$ with cardinality equal to $17,21,25,$ and $27$, respectively (see Figures~\ref{orderingfig} and~\ref{sresults1}). We deduce from the results that the interpolating polynomials corresponding to the incomplete sparse grids have stability and accuracy comparable to the respective complete sparse grids between which they lie. Moreover, we find that the interpolating polynomials over completed sparse grid filaments are the most stable and accurate ones of the lot.

\section*{Concluding remarks}

The application of the maximum volume principle in the selection of a well-defined and well-behaving polynomial basis for the Lagrange interpolation problem has been investigated in this paper. It has been demonstrated that the reciprocals of the volume as well as the minimum singular value of the Vandermonde matrix can be used to give an upper bound on the associated Lebesgue constant. In the framework of polynomial interpolation, the volume of the Vandermonde matrix thus has a natural interpretation as an indicator of the stability and accuracy of the interpolating polynomial, which has been observed in the numerical experiments conducted for both random and deterministic node configurations.

One could also consider fixing a set of low degree basis functions and computing the maximum volume submatrix of the generalized Vandermonde matrix with this predetermined set of basis functions removed. While this procedure may not result in the absolute (approximate) maximum volume Vandermonde submatrix, it restricts the number of high degree basis functions included in the interpolating polynomial---a desirable feature in practice---while still inheriting some stability endowed by the maximum volume approach. Such a procedure is not mutually exclusive to the results and methodology presented in this paper.

\section*{Acknowledgement}

The author is immensely grateful for the graceful advice and patience of Dr.~Harri Hakula and Prof.~Nuutti Hyv\"{o}nen throughout the long journey of writing this article.

\bibliographystyle{siam}
\bibliography{polyinterp}

\appendix

\section{Explicit orderings of the Smolyak node sequences and basis functions used in the numerical examples}\label{sappendix}

\begin{table}[!h]
\centering
\caption{Nodes and basis polynomials of the Smolyak--Clenshaw--Cutis rule with $d=2$ and $k=2$.}
\begin{tabular}{c|ll}
$i$&$\mathbf{x}_i^\textup{T}$&$\phi_i$\\ 
\hline 1& $[-1.000000, 0.000000]$ & $1$\\
2& $[-0.707107, 0.000000]$& $x_1$\\
3& $[0.000000, 0.000000]$& $-1 + 2 x_1^2$\\
4& $[0.707107, 0.000000]$& $x_2$\\
5& $[1.000000, 0.000000]$& $-1 + 2 x_2^2$\\
6& $[0.000000, -1.000000]$& $-3 x_1 + 4 x_1^3$\\
7& $[0.000000, -0.707107]$& $1 - 8 x_1^2 + 8 x_1^4$\\
8& $[0.000000, 0.707107]$& $-3 x_2 + 4 x_2^3$\\
9& $[0.000000, 1.000000]$& $1 - 8 x_2^2 + 8 x_2^4$\\
10& $[-1.000000, -1.000000]$& $x_1 x_2$\\
11& $[-1.000000, 1.000000]$& $x_1 (-1 + 2 x_2^2)$\\
12& $[1.000000, -1.000000]$& $(-1 + 2 x_1^2) x_2$\\
13& $[1.000000, 1.000000]$& $(-1 + 2 x_1^2) (-1 + 2 x_2^2)$
\end{tabular}
\end{table}

\begin{table}[!t]
\centering
\caption{Nodes and basis polynomials of the Smolyak--Clenshaw--Curtis rule with $d=2$ and $k=3$, where the nodes and basis functions corresponding to case $k=2$ have been removed.}
\begin{tabular}{c|ll}
$i$&$\mathbf{x}_i^\textup{T}$&$\phi_i$\\
\hline 14&$[-0.923880, 0.000000]$&$ 5 x_1 - 20 x_1^3 + 16 x_1^5$\\
 15&$[-0.382683, 0.000000]$&$ -1 + 18 x_1^2 - 48 x_1^4 + 32 x_1^6$\\
 16&$[0.382683, 
  0.000000]$&$ -7 x_1 + 56 x_1^3 - 112 x_1^5 + 64 x_1^7$\\
 17&$[0.923880, 0.000000]$&$ 
  1 - 32 x_1^2 + 160 x_1^4 - 256 x_1^6 + 128 x_1^8$\\
 18&$[0.000000, -0.923880]$&$ 5 x_2 - 20 x_2^3 + 16 x_2^5$\\
 19&$[0.000000, -0.382683]$&$ -1 + 18 x_2^2 - 48 x_2^4 + 32 x_2^6$\\
 20&$[0.000000, 0.382683]$&$ -7 x_2 + 56 x_2^3 - 112 x_2^5 + 
   64 x_2^7$\\
 21&$[0.000000, 0.923880]$&$ 
  1 - 32 x_2^2 + 160 x_2^4 - 256 x_2^6 + 128 x_2^8$\\
 22&$[-0.707107, -1.000000]$&$ (-3 x_1 + 4 x_1^3) x_2$\\
 23&$[-0.707107, 
  1.000000]$&$ (-3 x_1 + 4 x_1^3) (-1 + 2 x_2^2)$\\
 24&$[0.707107, -1.000000]$&$ (1 - 8 x_1^2 + 8 x_1^4) x_2$\\
 25&$[0.707107, 
  1.000000]$&$ (1 - 8 x_1^2 + 8 x_1^4) (-1 + 2 x_2^2)$\\
 26&$[-1.000000, -0.707107]$&$ x_1 (-3 x_2 + 4 x_2^3)$\\
 27&$[-1.000000, 0.707107]$&$ 
  x_1 (1 - 8 x_2^2 + 8 x_2^4)$\\
 28&$[1.000000, -0.707107]$&$ (-1 + 2 x_1^2) (-3 x_2 + 
     4 x_2^3)$\\
 29&$[1.000000, 
  0.707107]$&$ (-1 + 2 x_1^2) (1 - 8 x_2^2 + 8 x_2^4)$

\end{tabular}
\end{table}

\begin{table}[!t]
\centering
\caption{Nodes and basis polynomials of the Smolyak--Clenshaw--Cutis rule with $d=3$ and $k=2$.}
\begin{tabular}{c|ll}
$i$&$\mathbf{x}_i^\textup{T}$&$\phi_i$\\
\hline 1&$[-1.000000, 0.000000, 0.000000]$& $1$\\
 2&$[-0.707107, 0.000000, 0.000000]$& $x_1$\\
 3&$[0.000000, 0.000000, 0.000000]$& $-1 + 2 x_1^2$\\
 4&$[0.707107, 0.000000, 0.000000]$& $x_2$\\
 5&$[1.000000, 0.000000, 0.000000]$& $-1 + 2 x_2^2$\\
 6&$[0.000000, -1.000000, 0.000000]$& $x_3$\\
 7&$[0.000000, -0.707107, 0.000000]$& $-1 + 2 x_3^2$\\
 8&$[0.000000, 0.707107, 0.000000]$& $-3 x_1 + 4 x_1^3$\\
 9&$[0.000000, 1.000000, 0.000000]$& $1 - 8 x_1^2 + 8 x_1^4$\\
 10&$[0.000000, 0.000000, -1.000000]$& $-3 x_2 + 4 x_2^3$\\
 11&$[0.000000, 0.000000, -0.707107]$& $1 - 8 x_2^2 + 8 x_2^4$\\
 12&$[0.000000, 0.000000, 0.707107]$& $-3 x_3 + 4 x_3^3$\\
 13&$[0.000000, 0.000000, 1.000000]$& $1 - 8 x_3^2 + 8 x_3^4$\\
 14&$[-1.000000, -1.000000, 0.000000]$& $x_1 x_2$\\
 15&$[-1.000000, 1.000000, 0.000000]$& $x_1 (-1 + 2 x_2^2)$\\
 16&$[1.000000, -1.000000, 0.000000]$& $(-1 + 2 x_1^2) x_2$\\
 17&$[1.000000, 1.000000, 
  0.000000]$& $(-1 + 2 x_1^2) (-1 + 2 x_2^2)$\\
 18&$[-1.000000, 0.000000, -1.000000]$& $x_1 x_3$\\
 19&$[-1.000000, 0.000000, 1.000000]$& $x_1 (-1 + 2 x_3^2)$\\
 20&$[1.000000, 0.000000, -1.000000]$& $(-1 + 2 x_1^2) x_3$\\
 21&$[1.000000, 0.000000, 
  1.000000]$& $(-1 + 2 x_1^2) (-1 + 2 x_3^2)$\\
 22&$[0.000000, -1.000000, -1.000000]$& $x_2 x_3$\\
 23&$[0.000000, -1.000000, 1.000000]$& $x_2 (-1 + 2 x_3^2)$\\
 24&$[0.000000, 1.000000, -1.000000]$& $(-1 + 2 x_2^2) x_3$\\
 25&$[0.000000, 1.000000, 
  1.000000]$& $(-1 + 2 x_2^2) (-1 + 2 x_3^2)$
\end{tabular}
\end{table}

\begin{table}[!t]
\centering
\caption{Nodes and basis polynomials of the Smolyak--Clenshaw--Curtis rule with $d=3$ and $k=3$, where the nodes and basis functions corresponding to case $k=2$ have been removed.}
\begin{tabular}{c|ll}
$i$&$\mathbf{x}_i^\textup{T}$&$\phi_i$\\
\hline 26&$[-0.923880,0.000000, 0.000000]$&$ 5 x_1 - 20 x_1^3 + 16 x_1^5$\\
 27&$[-0.382683,0.000000, 0.000000]$&$ -1 + 18 x_1^2 - 48 x_1^4 + 32 x_1^6$\\
 28&$[0.382683,0.000000, 
  0.000000]$&$ -7 x_1 + 56 x_1^3 - 112 x_1^5 + 64 x_1^7$\\
 29&$[0.923880,0.000000, 0.000000]$&$ 
  1 - 32 x_1^2 + 160 x_1^4 - 256 x_1^6 + 128 x_1^8$\\
 30&$[0.000000, -0.923880, 0.000000]$&$ 5 x_2 - 20 x_2^3 + 16 x_2^5$\\
 31&$[0.000000, -0.382683, 0.000000]$&$ -1 + 18 x_2^2 - 48 x_2^4 + 32 x_2^6$\\
 32&$[0.000000, 0.382683, 
  0.000000]$&$ -7 x_2 + 56 x_2^3 - 112 x_2^5 + 64 x_2^7$\\
 33&$[0.000000, 0.923880, 0.000000]$&$ 
  1 - 32 x_2^2 + 160 x_2^4 - 256 x_2^6 + 128 x_2^8$\\
 34&$[0.000000,0.000000, -0.923880]$&$ 5 x_3 - 20 x_3^3 + 16 x_3^5$\\
 35&$[0.000000,0.000000, -0.382683]$&$ -1 + 18 x_3^2 - 48 x_3^4 + 32 x_3^6$\\
 36&$[0.000000,0.000000, 0.382683]$&$ -7 x_3 + 56 x_3^3 - 112 x_3^5 + 
   64 x_3^7$\\
 37&$[0.000000,0.000000, 0.923880]$&$ 
  1 - 32 x_3^2 + 160 x_3^4 - 256 x_3^6 + 128 x_3^8$\\
 38&$[-0.707107, -1.000000, 0.000000]$&$ (-3 x_1 + 4 x_1^3) x_2$\\
 39&$[-0.707107, 1.000000, 
  0.000000]$&$ (-3 x_1 + 4 x_1^3) (-1 + 2 x_2^2)$\\
 40&$[0.707107, -1.000000, 
  0.000000]$&$ (1 - 8 x_1^2 + 8 x_1^4) x_2$\\
 41&$[0.707107, 1.000000, 
  0.000000]$&$ (1 - 8 x_1^2 + 8 x_1^4) (-1 + 2 x_2^2)$\\
 42&$[-0.707107,0.000000, -1.000000]$&$ (-3 x_1 + 4 x_1^3) x_3$\\
 43&$[-0.707107,0.000000, 
  1.000000]$&$ (-3 x_1 + 4 x_1^3) (-1 + 2 x_3^2)$\\
 44&$[0.707107, 
  0.000000, -1.000000]$&$ (1 - 8 x_1^2 + 8 x_1^4) x_3$\\
 45&$[0.707107,0.000000, 
  1.000000]$&$ (1 - 8 x_1^2 + 8 x_1^4) (-1 + 2 x_3^2)$\\
 46&$[-1.000000, -0.707107, 0.000000]$&$ x_1 (-3 x_2 + 4 x_2^3)$\\
 47&$[-1.000000, 0.707107, 0.000000]$&$ 
  x_1 (1 - 8 x_2^2 + 8 x_2^4)$\\
 48&$[1.000000, -0.707107,
  0.000000]$&$ (-1 + 2 x_1^2) (-3 x_2 + 4 x_2^3)$\\
 49&$[1.000000, 0.707107, 
  0.000000]$&$ (-1 + 2 x_1^2) (1 - 8 x_2^2 + 8 x_2^4)$\\
 50&$[-1.000000,0.000000, -0.707107]$&$ x_1 (-3 x_3 + 4 x_3^3)$\\
 51&$[-1.000000,0.000000, 0.707107]$&$ 
  x_1 (1 - 8 x_3^2 + 8 x_3^4)$\\
 52&$[1.000000, 
  0.000000, -0.707107]$&$ (-1 + 2 x_1^2) (-3 x_3 + 4 x_3^3)$\\
 53&$[1.000000,0.000000, 
  0.707107]$&$ (-1 + 2 x_1^2) (1 - 8 x_3^2 + 8 x_3^4)$\\
 54&$[0.000000, -0.707107, -1.000000]$&$ (-3 x_2 + 4 x_2^3) x_3$\\
 55&$[0.000000, -0.707107, 
  1.000000]$&$ (-3 x_2 + 4 x_2^3) (-1 + 2 x_3^2)$\\
 56&$[0.000000, 0.707107, -1.000000]$&$ (1 - 8 x_2^2 + 
     8 x_2^4) x_3$\\
 57&$[0.000000, 0.707107, 
  1.000000]$&$ (1 - 8 x_2^2 + 8 x_2^4) (-1 + 2 x_3^2)$\\
 58&$[0.000000, -1.000000, -0.707107]$&$ x_2 (-3 x_3 + 4 x_3^3)$\\
 59&$[0.000000, -1.000000, 0.707107]$&$ 
  x_2 (1 - 8 x_3^2 + 8 x_3^4)$\\
 60&$[0.000000, 1.000000, -0.707107]$&$ (-1 + 2 x_2^2) (-3 x_3 + 
     4 x_3^3)$\\
 61&$[0.000000, 1.000000, 
  0.707107]$&$ (-1 + 2 x_2^2) (1 - 8 x_3^2 + 8 x_3^4)$\\
 62&$[-1.000000, -1.000000, -1.000000]$&$ 
  x_1 x_2 x_3$\\
 63&$[-1.000000, -1.000000, 1.000000]$&$ 
  x_1 x_2 (-1 + 2 x_3^2)$\\
 64&$[-1.000000, 1.000000, -1.000000]$&$ 
  x_1 (-1 + 2 x_2^2) x_3$\\
 65&$[-1.000000, 1.000000, 1.000000]$&$ 
  x_1 (-1 + 2 x_2^2) (-1 + 2 x_3^2)$\\
 66&$[1.000000, -1.000000, -1.000000]$&$ (-1 + 
     2 x_1^2) x_2 x_3$\\
 67&$[1.000000, -1.000000, 
  1.000000]$&$ (-1 + 2 x_1^2) x_2 (-1 + 2 x_3^2)$\\
 68&$[1.000000, 
  1.000000, -1.000000]$&$ (-1 + 2 x_1^2) (-1 + 
     2 x_2^2) x_3$\\
 69&$[1.000000, 1.000000, 
  1.000000]$&$ (-1 + 2 x_1^2) (-1 + 2 x_2^2) (-1 + 2 x_3^2)$

\end{tabular}
\end{table}

\end{document}